\documentclass[11pt, dina4]{article}
\usepackage{graphicx,epsfig}
\usepackage{epstopdf}
\usepackage{amssymb,amsmath}
\usepackage{cases}
\usepackage{amsthm}
\usepackage{caption,lipsum}
\usepackage{stmaryrd}
\usepackage{tabularx}
\usepackage{subfigure}
\usepackage{color, soul}
\usepackage{multirow}
\usepackage{empheq,float} 
\DeclareGraphicsExtensions{.eps}
\usepackage{hyperref}
\usepackage[margin=1in]{geometry} 
\usepackage{comment}

\newcommand{\p}{\partial}
\newcommand{\Og}{\Omega}
\newcommand{\fl}[2]{\frac{#1}{#2}}
\newcommand{\dt}{\delta}
\newcommand{\nn}{\nonumber}
\newcommand{\ap}{\alpha}
\newcommand{\bt}{\beta}

\newcommand{\veps}{\varepsilon}
\newcommand{\Dt}{\Delta}

\newcommand{\bea}{\begin{eqnarray}}
\newcommand{\eea}{\end{eqnarray}}
\newcommand{\beas}{\begin{eqnarray*}}
\newcommand{\eeas}{\end{eqnarray*}}
\newtheorem{remark}{Remark}[section]
\newtheorem{lemma}{Lemma}[section]

\newcommand{\bx}{{\bf x} }
\newcommand{\by}{{\bf y} }

\definecolor{ForestGreen}{rgb}{0.0, 0.5, 0.0}

\numberwithin{equation}{section}

\title{Variable-order fractional Laplacian and its accurate and efficient computations with meshfree methods}
\author{Yixuan Wu\thanks{Department of Pharmacology, University of California Davis, Davis, CA 95616-8636 (Email: ywxwu@ucdavis.edu)},  \ \
Yanzhi Zhang\thanks{Department of Mathematics and Statistics, Missouri University of Science and Technology, Rolla, MO 65409 (Email:  zhangyanz@mst.edu)}}

\begin{document}
\date{}
\maketitle

\begin{abstract}
The variable-order fractional Laplacian plays an important role in the study of heterogeneous systems. 
In this paper, we propose the first numerical methods for the variable-order Laplacian 
$(-\Dt)^{\ap(\bx)/2}$ with  $0 < \ap(\bx) \le 2$, 
which will also be referred as the variable-order fractional Laplacian if $\ap(\bx)$ is strictly less than $2$. 
We present a class of hypergeometric functions whose variable-order Laplacian can be analytically expressed. 
Building on these analytical results, we design the meshfree methods 
based on globally supported radial basis functions (RBFs), including Gaussian, generalized inverse multiquadric, and Bessel-type RBFs, to approximate the variable-order Laplacian   $(-\Dt)^{\ap(\bx)/2}$. 
Our meshfree methods integrate the advantages of both pseudo-differential and hypersingular integral forms of the variable-order fractional Laplacian, and thus avoid numerically approximating the hypersingular integral. 
Moreover, our methods are simple and flexible of domain geometry, and their computer implementation remains the same for any dimension $d \ge 1$. 
Compared to finite difference methods, our methods can achieve a desired accuracy with much fewer points. 
This fact makes our method much attractive for problems involving variable-order fractional Laplacian where the number of points required is a critical cost.  
We then apply our method to study solution behaviors of variable-order fractional PDEs arising in different fields, including transition of waves between classical and fractional media, and coexistence of anomalous and normal diffusion in both diffusion equation and the Allen–Cahn equation. 
These results would provide insights for further understanding and applications of variable-order fractional derivatives. 
\end{abstract}




\section{Introduction}
\label{section1}
\setcounter{equation}{0}

Recently,  variable-order fractional differential equations have attracted great attention in modeling heterogeneous properties of complex systems  \cite{Sun2019,  Lorenzo2002, Xiang2019, Zhu2014, Zhang2012, Xue2018, Delia2021}. 
In contrast to traditional (constant-order) derivatives,  the order of variable-order fractional derivatives may depend on time, space, or even dependent variables, which enables to easily study heterogeneous temporal or spatial effects. 
For example, coexistence and transition between anomalous and normal diffusion have been studied in many fields, including biology \cite{Javanainen2013}, turbulence \cite{Dubrulle1998}, and geophysics \cite{Baeumer2010, Meerschaert2008}, where variable-order fractional derivatives play an important role in describing such heterogeneous behaviors \cite{Lenzi2016, Dwivedi2021}. 
The variable-order fractional Laplacian $(-\Dt)^{\ap(\bx)/2}$ represents the infinitesimal generator of a stable-like Feller process, where  exponent $\ap(\bx)$ is usually spatial dependent  \cite{Kuhn2020, Chen2020, Kuhn2017, Samko1993}. 
Compared to its constant-order counterpart,  current studies on variable-order fractional Laplacian still remain  limited. 
Specifically,  a large number of functions have been reported in the literature (e.g. \cite{Dyda2012, Dyda2017} and references therein),  for which the action of constant-order fractional Laplacian can be written analytically. 
However,  no corresponding studies can be found on the variable-order fractional Laplacian. 
Moreover, the spatial-dependent exponent $\ap(\bx)$ makes numerical simulations extremely challenging, and no numerical method has been reported for the variable-order fractional Laplacian. 

The main purpose of this work is to fill these fundamental gaps in the literature. 
To this end, we will introduce the first numerical methods for the variable-order Laplacian 
and present a class of hypergeometric functions that their variable-order  Laplacian can be analytically expressed.
Let $\Og \subset {\mathbb R}^d$ (for $d \ge 1$) be an open bounded domain.  
Consider the variable-order fractional Poisson equation as follows \cite{Kuhn2017, Kuhn2020}:
\begin{equation}
\begin{split}
\label{diffusion}
(-\Dt)^{{\ap(\bx)}/{2}}u(\bx) = f(\bx), \qquad &\mbox{for} \ \, \bx \in \Og, \\ 
u(\bx)  = g(\bx), \qquad &\mbox{for} \ \, \bx \in \Og^c, 
\end{split}
\end{equation}
where $\Og^c = {\mathbb R}^d\backslash\Og$ represents the complement of domain $\Og$. 
The operator $(-\Dt)^{\ap(\bx)/2}$ is the variable-order Laplacian with spatial-dependent exponent $0 < \ap(\bx) \le 2$. 
For notational convenience, we will refer $(-\Dt)^{\ap(\bx)/2}$ as the {\it variable-order Laplacian} if $0 < \ap(\bx) \le 2$, or as the {\it variable-order  fractional Laplacian} if $0 < \ap(\bx) < 2$ (i.e., $\ap(\bx)$ is strictly less than $2$). 
More detailed discussion of this operator can be found in Section \ref{section2}. 

Generally, fractional derivatives of variable-order can be viewed as a heterogeneous generalization of their constant-order counterparts, such that the order can vary  as a function of dependent or independent variables. 
Since the order may depend on time, space, or even an independent external variables, the potential definitions of variable-order fractional derivatives could be vast. 
But, one can roughly classify them into four main types based on the definitions of their constant-order counterparts, including the Riesz derivatives, Gr\"{u}nwald--Letnikov derivatives, Riemann--Liouville derivatives, and Caputo derivatives \cite{Sun2019, Meerschaert2006,Lin2009,Zhao2015}. 
So far, many analytical results of variable-order operators can be found in the field of {\it variable exponent analysis}, where the parameters of operators and/or spaces may vary from point to point (instead of constant everywhere). 
However, compared to the constant-order derivatives, the current understanding of variable-order fractional derivatives still remains scant, and most of existing studies focus on the variable-order Riemann--Liouville derivatives \cite{Sun2019, Lin2009} and Caputo derivatives \cite{Zhao2015, Sun2019}. 

In this work, we propose the first numerical methods for the  variable-order Laplacian $(-\Dt)^{\ap(\bx)/2}$ with $0 \le \ap(\bx) \le 2$ and apply them to study solution behaviors of variable-order fractional PDEs. 
The main contributions of this work can be summarized as follows. 
\begin{itemize}\itemsep -1pt
\item[(i)]  We present a class of hypergeometric functions whose variable-order Laplacian can be analytically expressed.  
Noticing the relation between hypergeometric functions and many other functions, we further obtain analytical results for the Gaussian functions, generalized inverse multiquadric functions,  Bessel-type functions, and compactly supported functions on unit ball $B_1({\bf 0})$. 
This is the first time that such special functions are reported for the variable-order Laplacian. 
It is pointed out in \cite[p. 267]{Bonito2019} that ``{\it One of the difficulties in developing numerical approximation to (5) [fractional Poisson equation] is that there are relatively few examples where analytical solutions are available}". 
Our analytical results play an important role in studying the properties of variable-order Laplacian and also provide researchers a rich list of benchmark results for testing numerical methods. 
Furthermore, the analytical results on Gaussian, generalized inverse multiquadric, and Bessel-type functions  provide the foundation in the development of our meshfree radial basis function 
(RBF) methods. 

\item[(ii)] We propose the first numerical methods for the variable-order Laplacian  $(-\Dt)^{\ap(\bx)/2}$. 
More precisely,   three RBF-based numerical methods are introduced for this heterogeneous operator, including Gaussian RBFs, generalized inverse multiquadric RBFs, and Bessel-type RBF methods. 
Our meshfree methods combine the advantages of pseudo-differential  representation and hypersingular integral form of the variable-order fractional Laplacian and  thus bypass approximating the hypersingular integral in the fractional Laplacian. 
Consequently, they avoid large computational cost in evaluating the fractional derivative of RBFs, which is one main novelty distinguishing our method from other RBF-based methods \cite{Pang2015, Piret2013, Rosenfeld2019}. 
Moreover, our methods are simple and flexible of domain geometry, and their computer implementation remains the same for any dimension $d \ge 1$. 

\item[iii)] We numerically study the solution behaviors of various PDEs arising in modeling heterogeneous media.  
It is well known that the combination of nonlocality and heterogeneity introduce formidable  challenges in studying such systems. 
Numerical studies show that our methods are very effective in solving problems with  variable-order Laplacian. 
The transition of waves between classical and fractional media is studied in wave equations, while the coexistence of anomalous and normal diffusion are explored with the diffusion equation as well as the  Allen--Cahn equations.  
These studies could provide insights for the further understanding and applications of variable-order fractional derivatives. 
\end{itemize}

The paper is organized as follows. 
In Section \ref{section2},  two definitions of the variable-order fractional Laplacian are introduced together with some fundamental properties. 
Moreover, we present a collection of special functions for which the action of  variable-order Laplacian can be analytically expressed. 
In Section \ref{section3}, we propose a class of meshfree methods based on globally supported RBFs to discretize the variable-order Laplacian $(-\Dt)^{\ap(\bx)/2}$. 
The performance of our methods in approximating the Laplacian operators is tested in Section \ref{section4}, and various PDEs in modeling heterogeneous media are explored in Section \ref{section5}. 
Finally, we summarize the paper in Section \ref{section6}.

\section{Variable-order (fractional) Laplacian}
\label{section2}
\setcounter{equation}{0}

The variable-order fractional Laplacian  can be viewed as a heterogeneous generalization of the celebrated fractional Laplacian $(-\Dt)^{\fl{\ap}{2}}$ with  constant exponent $\ap \in (0, 2)$.  
It plays an important role in the study of heterogeneous problems. 
However, compared to the constant-order operator  $(-\Dt)^{\fl{\ap}{2}}$, the current understanding of  variable-order fractional Laplacian still falls very behind. 
In this section, we will first introduce the variable-order fractional Laplacian $(-\Dt)^{\ap(\bx)/2}$ in both pointwise integral form and pseudo-differential form,  and then discuss its properties in comparison to the constant-order Laplacian $(-\Dt)^{\fl{\ap}{2}}$.  
In Section \ref{section2-1}, we will present a collection of functions whose variable-order Laplacian can be analytically written. 
These analytical results can not only advance our understanding on the variable-order fractional Laplacian but also provide the key foundation in developing our numerical methods. 

The variable-order fractional Laplacian $(-\Dt)^{\ap(\bx)/2}$,  representing the infinitesimal generator of a  stable-like Feller process, can be defined in a hypersingular integral form \cite{Bass1988, Samko1993, Chen2020, Luo2019, Kuhn2020, Kuhn2017}:
\bea\label{integralFL}
(-\Delta)^{{\alpha(\bx)}/{2}}u(\bx) = C_{d,\alpha(\bx)}\,{\rm P. V.}\int_{{\mathbb R}^d}\fl{u(\bx) - u({\bf y})}{|\bx -{\bf y}|^{d+\ap(\bx)}}d{\bf y}, 
\eea
for $0 < \ap(\bx) < 2$, where ${\rm P. V.}$ stands for the principal value integral, and the normalization function is given by
\beas\label{constant}
C_{d,\alpha(\bx)}=\fl{2^{\ap(\bx)-1} \ap(\bx)\,\Gamma\big(\fl{{\ap(\bx)+d}}{2}\big)}{\sqrt{\pi^{d}}\,\Gamma\big(1 -\fl{\ap(\bx)}{2}\big)}
\eeas
with $\Gamma(\cdot)$ being the Gamma function.  
It is pointed out in \cite{Bass1988, Kuhn2017, Kuhn2020} that a sufficient condition for the variable-order fractional Laplacian being the infinitesimal generator of a stable-like process is that  $\ap: 
{\mathbb R}^d \to (0, 2)$ is H\"older continuous and $\inf_{\bx \in {\mathbb R}^d} \ap(\bx) > 0$. 
The integral in  (\ref{integralFL}) provides a pointwise definition of the variable-order fractional Laplacian, which can be viewed as an immediate generalization of the well-known (constant-order) fractional Laplacian $(-\Dt)^\fl{\ap}{2}$. 
It can be also obtained via the inverse of variable-order Riesz potential \cite{Samko1993, Samko2013}. 
The variable-order fractional Laplacian is a nonlocal operator where each point $\bx$ interacts with all the other points $\by \in {\mathbb R}^d$, but the kernel function characterizing interactions may vary  point to point (i.e., depending on $\ap(\bx)$). 
Note that the integral definition in (\ref{integralFL}) is valid for $0 < \ap(\bx) < 2$ and not compatible with the classical Laplacian operator $-\Dt = -\big(\p_{x^{(1)}}^2 + \p_{x^{(2)}}^2 + \cdots + \p_{x^{(d)}}^2\big)$. 

On the other hand,  the variable-order fractional Laplacian  is  often studied under the general framework of pseudo-differential operators.  
It can be defined as a pseudo-differential operator of symbol $|\xi|^{\ap(\bx)}$  \cite{Bass1988, Kikuchi1997, Leopold1999, Sun2019, SAMKO1993b, Kuhn2020}, i.e., 
\begin{eqnarray}
\label{pseudo}
(-\Delta)^{\alpha(\bx)/2}u({\bx}) = \int_{{\mathbb R}^d}\widehat{u}(\boldsymbol\xi)|\boldsymbol\xi|^{\alpha(\bx)} e^{2\pi i \bx\cdot\boldsymbol\xi} d\boldsymbol\xi, \qquad \mbox{for} \ \  \ap(\bx) > 0,
\end{eqnarray}
where $\widehat{u}(\boldsymbol\xi)$ represents the Fourier transform of $u(\bx)$. 
The pseudo-differential operator \eqref{pseudo} covers a wide class of Laplace operators. 
In the special case of $\ap(\bx) \equiv 2$,  it reduces to the spectral representation of the classical { negative} Laplacian $-\Dt$. 
If $\ap(\bx) \equiv \ap$ is a constant and $\ap \in (0, 2)$, it collapses to the well-known (constant-order) fractional Laplacian $(-\Dt)^\fl{\ap}{2}$. 
Moreover, the $\ap(\bx)$-parametric pseudo-differential operator in (\ref{pseudo}) unifies the classical and fractional Laplacians in a seamless way, enabling it to naturally describe the coexistence of normal ($\ap = 2$) and anomalous ($\ap < 2$) diffusion phenomena. 

In this work, we will focus on the variable-order Laplacian covering both classical and fractional Laplacians. 
As mentioned previously, we will refer $(-\Dt)^{\ap(\bx)/2}$ as the {\it variable-order fractional Laplacian } if  $0 < \ap(\bx) < 2$, or the {\it variable-order Laplacian }\,if  $0 < \ap(\bx) \le 2$. 
We assume that $0 < \inf  \ap(\bx) \le \ap(\bx) \le 2$ and $\ap(\bx)$ is H\"older continuous \cite{Bass1988, Kuhn2017, Kuhn2020}. 
More discussion on exponent $\ap(\bx)$ can be found in \cite{Bass1988, Samko1993, Samko2013, Luo2019, Chen2020} and references therein. 
It shows in \cite{Bass1988, SAMKO1993b, Kuhn2020, Kuhn2017} that the integral definition (\ref{integralFL}) and pseudo-differential definition  \eqref{pseudo} of the variable-order fractional Laplacian are equivalent for $0 < \ap(\bx) < 2$ and function $u \in C^\infty ({\mathbb R}^d)$.  
The pseudo-differential definition in (\ref{pseudo})  unifies the integer-order (i.e. { $\ap = 2m$ with $m \in {\mathbb N}$}) and fractional-order Laplacians in a single form via exponent $\ap(\bx)$. 
As we will see in Section \ref{section3}, this property plays a key role in developing compatible schemes for classical  and  fractional Laplacians, but it is challenging to incorporate non-periodic boundary conditions into the pseudo-differential form (\ref{pseudo}). 
In contrast, the integral definition in \eqref{integralFL} can easily work with non-periodic boundary conditions, but it is incompatible to the classical Laplacian, i.e. $\ap(\bx) \neq 2$ in \eqref{integralFL}. 
This motivates us to combine the advantages of both definitions such that we can study the variable-order Laplacian $(-\Dt)^{\ap(\bx)/2}$ but free of periodic boundary-condition constraints. 

The study of variable-order fractional Laplacian  can be traced back to the seminal paper \cite{Bass1988} where Bass studied pure jump Markov processes associated with such a generator. 
The existence and well-posedness of the martingale solution of $(-\Dt)^{\ap(\bx)/2}$ are studied in  \cite{Bass1988, Tsuchiya1992}, and it shows that there is a strong Markov process corresponding to the variable-order fractional Laplacian. 
Recently,  the Schauder estimates for the variable-order fractional Poisson equation are studied in \cite{Kuhn2020, Kuhn2017}, and solution properties of a more general elliptic problem with the integral operator
 (\ref{integralFL}) can be found in \cite{Xiang2019}.  
On the other side, the calculus of singular integral operators has been studied by recasting them into a theory of pseudo-differential operators \cite{Kohn1965, Hormander1965, Kikuchi1997, Leopold1999}. 
So far, numerous research can be found on {\it variable exponent analysis}, where operators are studied in variable exponent setting, i.e., the parameters defining the operators and/or the space may vary from point to point (see \cite{Rafeiro2016, CRUZ-URIBE2013, DIENING2017} and references therein). 
Some properties for the variable-order fractional Laplacian have been discussed under this framework. 

\subsection{Properties of variable-order Laplacian}
\label{section2-1}

Recently, numerous studies have been reported on the fractional derivatives with variable order \cite{Sun2019, Zhao2015, Song2016, Du2022}, including the Gr\"unwald--Letnikov derivatives, Riemann--Liouville derivatives, and Caputo derivatives. 
Compared to these fractional derivatives, analytical and numerical studies on the variable-order Laplacian still remain limited. 
Even though many results can be found on the constant-order fractional Laplacian $(-\Dt)^\fl{\ap}{2}$,  it is challenging to generalize them into the variable order cases.  
In fact the variable-order fractional Laplacian may lose some important properties (e.g., rotational invariance) of its constant-order counterpart. 

In the following, we will study the properties of variable-order Laplacian and present a collection of functions whose variable-order Laplacian can be analytically expressed. 
We remark that this section does not attempt to make a comprehensive study on the variable-order Laplacian, but instead concentrates on some important functions that $(-\Dt)^{\ap(\bx)/2}u$ can be found analytically. 
For notational convenience, let's denote 
\beas
{\mathcal U}(\bx, \by) := (-\Dt)^{\ap(\bx)/2}u(\by), \qquad \mbox{for} \ \ 0 < \ap(\bx) \le 2. 
\eeas
First, we present the following properties of variable-order Laplacian by generalizing those of the  constant-order fractional Laplacian in \cite{Burkardt2021}.
\begin{lemma}\label{lemma1} 
The variable-order Laplacian operator $(-\Dt)^{\ap(\bx)/2}$ satisfies the following properties: 
\bea
\label{prop1}
&&\qquad(-\Dt)^{\ap(\bx)/2}\big[u(\by - \by_0)\big] = {\mathcal U}(\bx,\, \by - \by_0),  \quad\ \ \mbox{for} \ \, \by_0\in{\mathbb R}^d, \\
\label{prop2}
&&\qquad(-\Dt)^{{\ap(\bx)}/{2}}\big[u(\zeta \by)\big] = |\zeta|^{\ap(\bx)}{\mathcal U}(\bx, \zeta \by), \quad\quad \ \,\mbox{for} \,\ \zeta\in{\mathbb R}.
\eea
\end{lemma}
The properties in (\ref{prop1})--(\ref{prop2})  are consistent with those of the classical ($\ap \equiv 2$) and constant-order fractional ($\ap < 2$) Laplacians. 
They play a fundamental role in the design of RBF-based numerical methods for variable-order Laplacian. 

Recently, the fractional Laplacian of Meijer G-functions and generalized hypergeometric functions have been extensively studied in \cite{Dyda2012, Dyda2017} for constant-order cases,  i.e., $\ap(\bx) \equiv {\ap}$ { and }$ { \ap } \in (0, 2)$.  
These results can be used to study the eigenvalues and eigenfunctions of the fractional Laplacian $(-
\Dt)^\fl{\ap}{2}$ in unit balls and also serve as benchmark results for testing numerical methods. 
In the following, we will generalize these results from the constant-order fractional Laplacian  $(-\Dt)^\fl{\ap}{2}$  for $\ap \in (0, 2)$ to the variable-order Laplacian $(-\Dt)^{\ap(\bx)/2}$ for $0 < \ap(\bx) \le 2$. 
We will mainly focus on  the generalized hypergeometric functions. 
Note that similar generalizations can be done for the Meijer G-functions under appropriate conditions, but we will leave them for  future study.

For convenience of discussion, we adopt the notations used in \cite{Dyda2017}, and assume $V(\bx)$ as a solid (homogeneous) harmonic polynomial of degree $l \in {\mathbb N}^0$, i.e., a polynomial satisfying $\Dt V = 0$ and  homogeneous of degree $l$, {with ${\mathbb N}^0$ denotes the set of nonnegative integers.
Let  $p, q \in {\mathbb N}^0$} and  $p \le q+1$.  Define the generalized hypergeometric function as 
\bea\label{hyper}
_pF_q\big((a_1, a_2, \cdots, a_p); \, (b_1, b_2, \cdots, b_{q}); \,r \Big) := \sum_{n=0}^{\infty} \fl{(a_1)_n \,(a_2)_n\, \cdots\, (a_p)_n}{(b_1)_n \,(b_2)_n\, \cdots\, (b_q)_n}\ \fl{r^n}{n!},\qquad 
\eea
where we assume that none of  $b_k$ is a nonpositive integer,  and the Pochhammer symbol {$(c)_{n} = c(c+1) \cdots (c+n-1)$ denotes the rising factorial with $(c)_0 = 1$.} 
If $p \le q$, the above series is convergent for all $r \in {\mathbb C}$. 
More information of hypergeometric functions can be found in \cite{GRADSHTEYN2007, PRUDNIKOV1990}. For the hypergeometric function $_pF_q$ and its variation, we have the following lemma:

\begin{lemma}[Laplacian of generalized hypergeometric functions]
\label{lemma2} 
Let  $p, q \in {\mathbb N}^0$ and $q-1 \le p \le q+1$.
Denote ${\bf a} = (a_1, a_2, \cdots, a_p)$ and \,${\bf b} = (b_1, b_2, \cdots, b_{q-1})$, and define function
\bea\label{hyper-u}
u(\bx) = V(\bx)\,_pF_q\Big({\bf a}; \, {\bf b}, \,\vartheta; \, -|\bx|^2\Big), \quad \ \ \mbox{for} \ \,\, \bx \in {\mathbb R}^d,
\eea
where $V(\bx)$ is a solid harmonic polynomial of degree { $l \in {\mathbb N}^0$}, and $\vartheta = \fl{d}{2}+l$. 
Then the variable-order Laplacian of $u(\bx)$ can be analytically given by
\bea\label{hyperu2}
(-\Dt)^{\ap(\bx)/2}u(\bx) = C_1(\bx)V(\bx)\,_pF_q\Big({\bf a} + \fl{\ap(\bx)}{2}; \, {\bf b} + \fl{\ap(\bx)}{2}, \,\vartheta; \, -|\bx|^2\Big),\quad \
\eea
where the coefficient 
\beas
C_1(\bx) = 2^{\ap(\bx)} \prod_{k=1}^p \fl{\Gamma\big(a_k + \fl{\ap(\bx)}{2}\big)}{\Gamma(a_k)}\prod_{k = 1}^{q-1}\fl{\Gamma(b_k)}{\Gamma\big(b_k+\fl{\ap(\bx)}{2}\big)}.
\eeas
\end{lemma}
In Lemma \ref{lemma2}, we not only generalize the result of  constant-order fractional Laplacian in  \cite[Corollary 2]{Dyda2017} to variable-order exponent $\ap(\bx)$, but also extend it  to include $\ap(\bx) = 2$. 
The proof of Lemma \ref{lemma2} can be done by following similar arguments for the constant-order fractional Laplacian in \cite{Dyda2011, Dyda2012, Dyda2017}. 
Noticing that many elementary functions can be written in terms of  hypergeometric functions, thus we can further obtain the following results.
\begin{lemma}[Laplacian of globally supported functions]
\label{lemma3} 
Let $V(\bx)$ be a solid harmonic polynomial of degree { $l \in {\mathbb N}^0$}. 
Denote $\vartheta = \fl{d}{2}+l$. 
The variable-order Laplacian of the following infinitely differentiable functions can be analytically expressed in terms of  hypergeometric functions. \vspace{1mm}

\begin{itemize}\itemsep 3pt
\item[(i)] For the Gaussian type function $u(\bx)= V(\bx) e^{-|\bx|^2}$, it holds
\bea\label{Gaussian}
(-\Dt)^{\ap(\bx)/2}u(\bx) = \frac{2^{\ap(\bx)}\Gamma\big(\vartheta + \fl{\ap(\bx)}{2}\big)}{\Gamma\big(\vartheta\big)}V(\bx)  \,_1F_1\Big(\vartheta + \fl{\ap(\bx)}{2}; \,\vartheta;\, -|\bx|^2\Big),\qquad
\eea
where $_1F_1$ represents the confluent hypergeometric function.   
\item[(ii)] For the inverse multiquadric type function $u(\bx) =  V(\bx)(1+|\bx|^2)^{-\bt}$ with $\bt > 0$, it holds
\bea\label{IMQ}
&&(-\Dt)^{\ap(\bx)/2}u(\bx) = \frac{2^{\ap(\bx)}\Gamma\big(\vartheta + \fl{\ap(\bx)}{2}\big)\Gamma\big(\bt + \fl{\ap(\bx)}{2}\big)}{\Gamma\big(\vartheta\big)\Gamma(\bt)}\cdot \nn\\
&&\hspace{3.5cm}V(\bx) \,_2F_1\Big(\vartheta + \fl{\ap(\bx)}{2}, \bt+\fl{\ap(\bx)}{2}; \, \vartheta; \, -|\bx|^2\Big),\qquad\qquad
\eea
where $\,_2F_1$ denotes the Gauss hypergeometric function. 
\vspace{1mm}
 
\item[(iii)]  For function $u(\bx) = V(\bx) {J_{s-1}(|\bx|)}/{|\bx|^{s-1}}$ with $J_{s-1}$  the Bessel function of order $(s -1)$ for $s > 0$, it holds
\bea\label{Bessel}
&&(-\Dt)^{\ap(\bx)/2}u(\bx)
 = \fl{2^{1-s}\Gamma\big(\vartheta + \fl{\ap(\bx)}{2}\big)}{\Gamma\big(\vartheta\big)\Gamma\big(s+\fl{\ap(\bx)}{2}\big)}\cdot \nn\\
 &&\hspace{3.5cm}V(\bx)\,_1F_2\Big(\vartheta + \fl{\ap(\bx)}{2}; \, s+\fl{\ap(\bx)}{2}, \vartheta; \, -\fl{1}{4}|\bx|^2\Big). \qquad\qquad
\eea
\end{itemize}
\end{lemma}

\begin{figure}[htb!]
\centerline{
\includegraphics[height = 4.6cm, width = 6.20cm]{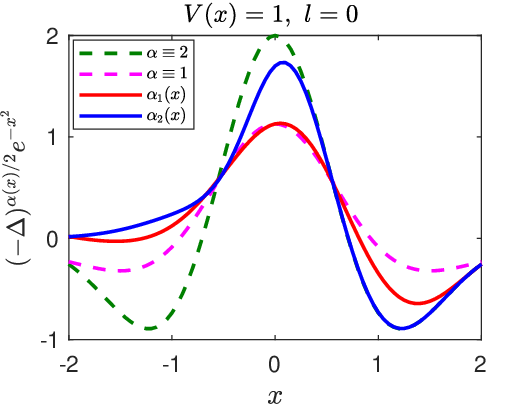}\hspace{-2mm}
\includegraphics[height = 4.6cm, width = 6.20cm]{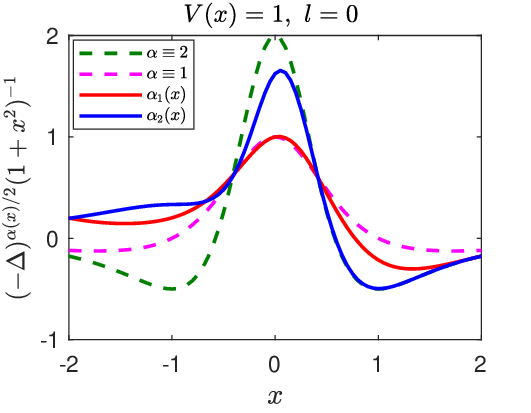}}
\centerline{
\includegraphics[height = 4.6cm, width = 6.20cm]{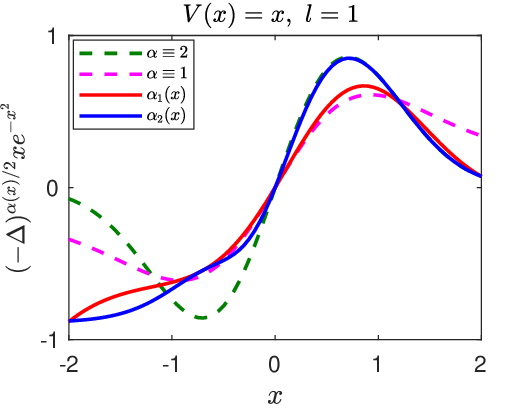}\hspace{-2mm}
\includegraphics[height = 4.6cm, width = 6.20cm]{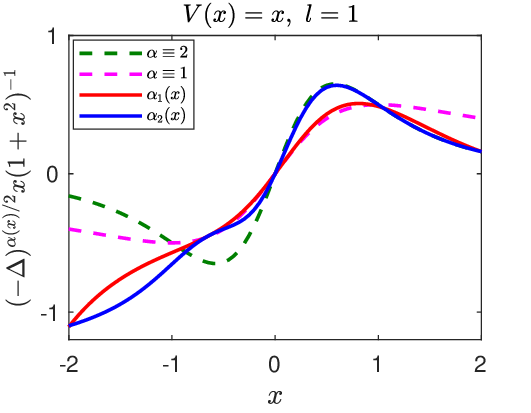}}
\caption{Illustration of variable-order Laplacian of the one-dimensional Gaussian type function (left column) and generalized inverse multiquadric type function with $\bt = 1$ (right column) for { $x \in (-2, 2)$}, where $\ap_1(x) = 1+x/2$ and $\ap_2(x) = 1+{\rm tanh}(2x+1)$.}\label{Figure2-1} 
\end{figure}

The results in Lemma \ref{lemma3} can be obtained from Lemma \ref{lemma2} by first rewriting the functions in (i)--(iii) in terms of  hypergeometric function $\,_pF_q$ and then applying the following relation:
\beas
\,_pF_q\big({\bf a}; \, {\bf b}; \, r\big) =  \,_{p+1}F_{q+1}\big({\bf a}, \, \eta; \, {\bf b},\,  \eta; \, r\big) = \,_{p+1}F_{q+1}\big(\eta, \, {\bf a}; \, \eta, \, {\bf b}; \, r\big),\quad \mbox{for} \ \ \eta \ge 0. 
\eeas
Specifically, we can write the Gaussian function, inverse multiquadric function, and Bessel-type functions as
\beas
e^{-|\bx|^2} = \,_0F_0\big(\ ;\ ; -|\bx|^2\big), && \\
\big(1 + |\bx|^2\big)^{-\bt}= \,_1F_0\Big(\bt; \,  ;  \, -|\bx|^2\Big), &&\mbox{for} \ \ \bt > 0,\\
\fl{J_{s-1}(|\bx|)}{|\bx|^{s-1}} = \frac{2^{1-s}}{\Gamma(s)}\,_0F_1\Big(\, ;\, s;\, -\fl{1}{4}|\bx|^2\Big), &\qquad & \mbox{for} \ \ s > 0.
\eeas
Figure \ref{Figure2-1} illustrates the heterogeneous properties of variable-order Laplacian acting on the one-dimensional Gaussian function and generalized inverse multiquadric function with $\bt = 1$. 
It shows that the results from constant order $\ap$ are symmetric (resp. antisymmetric) about $x = 0$ for $l = 0$ (resp. $l = 1$), owing to the rotational invariance  of  constant-order Laplacian. 
In contrast, the results from variable-order Laplacian may  lose this symmetry, depending on exponent $\ap(x)$. 

To the best of our knowledge, the results in Lemmas \ref{lemma2}--\ref{lemma3} are the first report on functions whose variable-order Laplacian can be analytically written. 
These analytical results can not only advance the understanding of  variable-order Laplacian but also serve as benchmarks in testing numerical methods for this operator. 
Note that lacking of benchmark results is one main challenge in the current literature  \cite{Bonito2019}. 
On the other hand, the Gaussian function (i.e., $e^{-r^2}$), inverse multiquadric functions (i.e., $(1+r^2)^{-\bt}$), and Bessel-type functions (i.e., $J_{s-1}(r)/r^{s-1}$) are well studied in the field of radial basis functions. 
They are the few positive definite functions among all RBFs. 
Hence, the results in Lemma \ref{lemma3} play a foundational role in the design of RBF-based meshfree methods for  variable-order Laplacian; see more discussion in Section \ref{section3}. 

If $s = (2m+1)/2$ with $m \in {\mathbb N}^0$, the function $J_{s-1}(|\bx|)/|\bx|^{s-1}$ can be alternatively expressed by means of regular trigonometric functions. 
Hence, we can further find the variable-order Laplacian of some important trigonometric functions from the results in (\ref{Bessel}). For example,  if $s = {1}/{2}$,  it holds $J_{-\fl{1}{2}}(|\bx|)\sqrt{|\bx|} = \sqrt{2/\pi}\cos(|\bx|)$, and thus 
\beas
(-\Dt)^{\ap(\bx)/{2}}\cos(|\bx|) = \fl{\sqrt{\pi} \Gamma\big(\fl{d + \ap(\bx)}{2}\big)}{\Gamma\big(\fl{d}{2}\big)\Gamma\big(\fl{1+\ap(\bx)}{2}\big)}\,_1F_2\Big(\fl{d+\ap(\bx)}{2}; \,\fl{1+\ap(\bx)}{2}, \fl{d}{2}; \,-\fl{1}{4}|\bx|^2\Big),
\eeas
It immediately implies that in one-dimensional (i.e., $d = 1$) case, there is 
\beas
(-\Dt)^{\ap(x)/{2}}\cos\big(x\big) = \cos\big(x\big), \qquad \mbox{for} \ \ x \in {\mathbb R}.
\eeas
If $s = {3}/{2}$, it holds $J_{\fl{1}{2}}(|\bx|)/\sqrt{|\bx|} = \sqrt{2/{\pi}}\sin(|\bx|)/|\bx|$, and thus
\beas
(-\Dt)^{\ap(\bx)/{2}}{\rm sinc}(|\bx|)= \fl{\pi^{\ap(\bx)+\fl{1}{2}}\Gamma\big(\fl{d+\ap(\bx)}{2}\big)}{2\Gamma\big(\fl{d}{2}\big)\Gamma\big(\fl{3+\ap(\bx)}{2}\big)}\,_1F_2\Big(\fl{d+\ap(\bx)}{2}; \, \fl{3+\ap(\bx)}{2}, \fl{d}{2}; \, -\fl{\pi^2}{4}|\bx|^2\Big), 
\eeas
where we define ${\rm sinc}(x) = \sin(\pi x)/(\pi x)$. 
One can continue and obtain analytical results for other $s = (2m+1)/2$, which we will omit here for the purpose of conciseness. 

In \cite{Dyda2012},  two compactly supported functions on a unit ball $B_1({\bf 0})$ have been analytically studied for the constant-order fractional Laplacian $(-\Dt)^\fl{\ap}{2}$. 
It shows in \cite{Dyda2012, Duo2015} that these results play an important role in studying the eigenvalues and eigenfunctions of the fractional Laplacian $(-\Dt)^\fl{\ap}{2}$. 
In Lemma \ref{lemma4}, we will generalize these results to the variable-order Laplacian. 
Notice that 
\beas
\big(1- r\big)^p = \,_1F_0 \big(-p; \, ; \, r\big), \qquad\mbox{for} \ \ |r|<1. 
\eeas
Then, we obtain the variable-order Laplacian of the following compactly supported functions.
\begin{lemma}[Laplacian of compactly supported functions  on $B_1({\bf 0})$]  \label{lemma4} 
Let $V(\bx)$ be a solid harmonic polynomial of degree { $l \in {\mathbb N}^0$}.  
Denote $\vartheta = \fl{d}{2}+l$. 
Suppose $u(\bx) = V(\bx)\big[(1-\lvert  \bx \rvert^2)^p\big]_+$ with constant $p > -1$. 
Then the variable-order Laplacian of $u(\bx)$ is given by
\bea\label{comp1}
&&(-\Dt)^{\ap(\bx)/2}u(\bx) = \fl{2^{\ap(\bx)}\Gamma(p+1)\Gamma\big(\vartheta + \fl{\ap(\bx)}{2}\big)}{\Gamma\big(\vartheta\big)\,\Gamma\big(p+1-\fl{\ap(\bx)}{2}\big)}\cdot\nn\\
&&\hspace{3.5cm}V(\bx) \,_2F_1 \Big(\vartheta + \fl{\ap(\bx)}{2},\, -p+\fl{\ap(\bx)}{2}; \, \vartheta; \, |\bx|^2\Big),\qquad\qquad 
\eea
for $|\bx| < 1$.
\end{lemma}
The proof of Lemma \ref{lemma4} can be done by following the similar lines in \cite{Dyda2012, Dyda2017}. 
The compactly supported functions in Lemma \ref{lemma4} are defined in a unit ball $B_1({\bf 0})$, but our result in (\ref{comp1}) can be generalized  to any ball $B_r({\bx_c})$ with $r > 0$ and $\bx_c\in {\mathbb R}^d$ by using the properties (\ref{prop1})--(\ref{prop2})  in Lemma \ref{lemma1}.  
Note that the results in \cite[Theorem 1]{Dyda2012} can be viewed as special cases of Lemma \ref{lemma4} by setting constant exponent $\ap \in (0, 2)$ and choosing different $V(\bx)$.   
Specifically, choosing $V(\bx) = 1$ leads to function $u(\bx) = \big[(1-|\bx|^2)^p\big]_+$, while  $V(\bx) = x^{(i)}$ gives $u(\bx) = x^{(i)}\big[(1-|\bx|^2)^p\big]_+$ with $\bx =\big (x^{(1)}, x^{(2)}, \cdots, x^{(d)}\big)$. 
\begin{figure}[htb!]
\centerline{
\includegraphics[height = 4.6cm, width = 6.20cm]{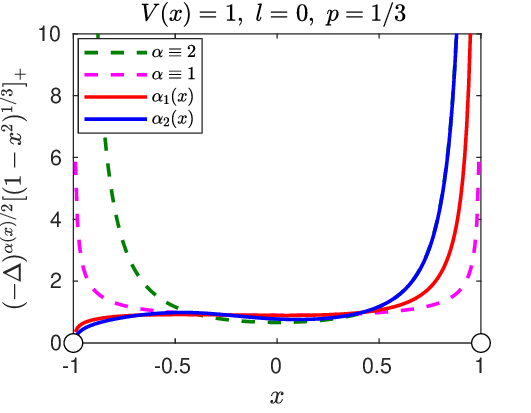}\hspace{-1mm}
\includegraphics[height = 4.6cm, width = 6.20cm]{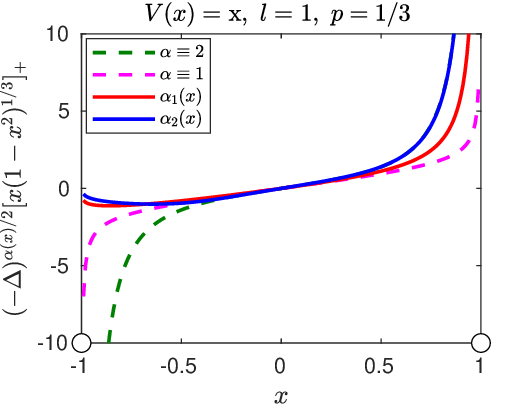}}
\centerline{
\includegraphics[height = 4.6cm, width = 6.20cm]{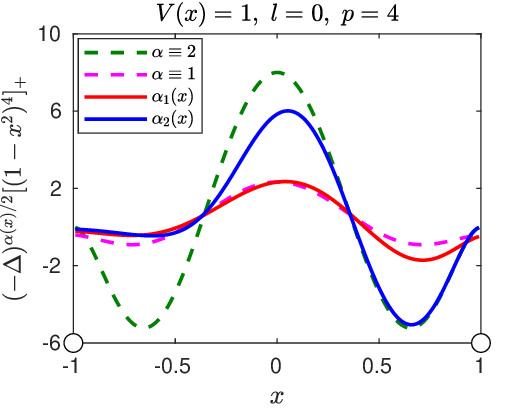}\hspace{-1mm}
\includegraphics[height = 4.6cm, width = 6.20cm]{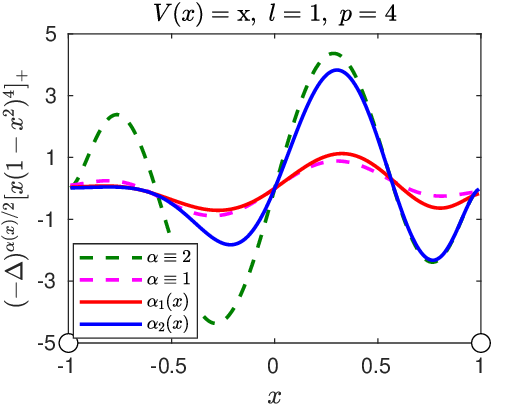}}
\caption{ Illustration of variable-order Laplacian of the one-dimensional compactly supported functions in Lemma \ref{lemma4}, where we choose $\ap_1(x) = 1+x/2$ and $\ap_2(x) = 1+{\rm tanh}(2x+1)$. The symbol `{\large$\circ$}' indicates that points $x = \pm 1$ are not included in the above plots.} \label{Figure2-2}
\end{figure}
Figure \ref{Figure2-2} shows the variable-order Laplacian of one-dimensional compactly supported functions for different $V(x)$ and $p$. 
Similar to our observations in Figure \ref{Figure2-1}, the results from constant order $\ap$ are symmetric (resp. antisymmetric) about $x = 0$ for $l = 0$ (resp. $l = 1$). 
Different from those in Lemma \ref{lemma3}, the compactly support functions here have finite smoothness at boundary (e.g., $x = \pm 1$), depending on the value of $p$. 
Hence, their Laplacian might go to infinity at boundary if $p$ is small (see, e.g., Figure \ref{Figure2-2} for $p = 1/3$).
\begin{remark} 
The variable-order Laplacian $(-\Dt)^{\ap(\bx)/2}$ covers a wide class of Laplace operators, including classical { negative} Laplacian $-\Dt$ (i.e. $\ap \equiv 2$) and constant-order fractional Laplacian $(-\Dt)^\fl{\ap}{2}$ for $\ap \in (0, 2)$. 
Hence, the analytical results discussed in Lemmas \ref{lemma2}--\ref{lemma4} hold for both classical and constant-order fractional Laplacians.  
\end{remark}

The recent literature \cite{Dyda2012, Dyda2017} provides an extensive list of special functions whose constant-order fractional Laplacian $(-\Dt)^\fl{\ap}{2}$ for $\ap \in (0, 2)$ can be analytically found.  
Most of their results can be generalized to the variable-order Laplacian with $0 < \ap(\bx) \le 2$, if the function $u$ is $\ap$-independent. 
However,  the situation becomes more complicated when $u$ is $\ap$-dependent, and  generalization from constant order to variable order in this case may not be straightforward. 
Further studies will be carried out in the future. 

\section{Meshfree RBF methods}
\label{section3}
\setcounter{equation}{0}

In this section, we will introduce our meshfree  methods based on radial basis functions (RBFs) to discretize the variable-order Laplacian $(-\Dt)^{\ap(\bx)/2}$. 
So far, many numerical methods have been proposed to approximate the constant-order fractional Laplacian $(-\Dt)^\fl{\ap}{2}$; see \cite{Acosta2017,  Ainsworth2018B, Duo2018, Duo-FDM2019, Duo-TFL2019, Bonito2019, Rosenfeld2019, Burkardt2021, Wu0020, Hao2021, Tang2020, Antoine2021} and references therein. 
In contrast, numerical methods for the variable-order Laplacian still remain very rare, and the main challenges stem from both its nonlocality and heterogeneity. 
In recent studies \cite{Zhu2014, Xue2018, Chen2014b}, the variable-order fractional Laplacian $(-\Dt)^{\ap(\bx)/2}$ are approximated by an averaged constant-order fractional Laplacian $(-\Dt)^{\fl{\bar{\ap}}{2}}$ with constant $\bar{\ap} = \int \ap(\bx) d\bx$, which makes numerical simulations much easier but loses the heterogeneous features of original models. 

Here, we propose our meshfree methods for the variable-order Laplacian $(-\Dt)^{\ap(\bx)/2}$, which enable us to effectively study heterogeneous phenomena in many applications \cite{Zhu2014, Xue2018, Chen2014b, Sun2019, Lenzi2016}. 
To the best of our knowledge, this is the first numerical methods developed for the variable-order fractional Laplacian $(-\Dt)^{\ap(\bx)/2}$. 
The key of our methods is to utilize the equivalence of pseudo-differential form and integral representation of the variable-order fractional Laplacian when acting on globally supported RBFs in (\ref{GA})--(\ref{BRBF}). 
Hence, we can combine the advantages of both definitions and bypass numerical approximation of the hypersingular integral in (\ref{integralFL}), which not only reduces computational cost but significantly simplifies the implementation  especially in high ($d > 1$) dimensions. 

We will start with a brief introduction of RBFs. 
RBFs are well-known for their success in function reconstruction from high-dimensional  scattered data.
They have been widely applied in many fields, including solving PDEs \cite{Kansa90I, Kansa90II}. 
RBFs are usually real-valued scalar functions that depend on the distance of point $\bx$ to a given center point $\bx_c$, i.e., $\varphi(\bx)=\varphi(|\bx- \bx_c|)$ for $\bx,\, \bx_c \in {\mathbb R}^d$. 
RBFs can be divided into two main categories:  globally supported functions and compactly supported functions. 
Among all globally supported RBFs, there is a class of infinitely differentiable positive definite functions, including 
\bea
\label{GA}
\mbox{Gaussian RBF:}   &\ &
\displaystyle \varphi(r) = e^{-r^2}, \\
\label{GIMQ}
\mbox{Generalized inverse multiquadric RBF:}  & & 
\displaystyle \varphi_\bt(r) = \big(1+r^2\big)^{-\bt}, \qquad \\ 
\label{BRBF}
\mbox{Bessel-based RBF:}  & & 
\displaystyle \varphi_m(r) = \fl{J_{m/2-1}(r)}{r^{m/2-1}}, 
\eea
where $\bt > 0$, $m \in {\mathbb N}$, and we denote $r = |\bx - \bx_c|$. 
Usually,  a shape parameter $\varepsilon > 0$ is introduced in globally supported RBFs,  and the functions are written as $\varphi^\varepsilon(r) := \varphi(\varepsilon r)$. 
The shape parameter $\varepsilon$ plays an important role in RBF-based methods, which is usually chosen as a constant.  
The variable shape parameters (e.g., depending on the center points) were also studied in the literature \cite{Fornberg2017, Wu0020}. 
Compared to the Gaussian RBFs, the studies of generalized inverse multiquadric (gIMQ) and Bessel-type RBFs (also referred to as oscillatory RBFs in \cite{Fornberg2006}) are still very recent. 
More discussion of RBFs can be found in \cite{Sarra2009, Flyer2006, Larsson2005} and references therein. 

Next, we will introduce our method to approximate the variable-order Laplacian $(-\Dt)^{\ap(\bx)/2}$ with extended Dirichlet boundary conditions in (\ref{diffusion}). 
Let $N$ and $\bar{N}$ be two positive integers, and $N < \bar{N}$. 
Denote $\bx_i \in \bar{\Og}$ (for $1 \le i \le {\bar{N}}$) as RBF center points on  $\bar{\Og} = \Og\cup\p\Og$. 
Specifically,  we let $\bx_i \in \Og$ for $1 \le i \le N$, and $\bx_i \in \p\Og$ for $N+1 \le i \le \bar{N}$. 
Assume that $u(\bx)$  can be approximated by
\bea\label{Sol1D}
\widehat{u}(\bx) :=\sum_{i=1}^{\bar{N}} \lambda_i\,\varphi^\veps(|\bx - \bx_i|),
\eea
where $\varphi^\veps(|\bx-\bx_i|)$ represents a RBF centered at point $\bx_i$. 
{ We consider the constant shape parameter $\veps$, i.e., all basis functions  $\varphi^\veps(|\bx-\bx_i|)$ in (\ref{Sol1D}) have the same  parameter $\veps$. 
Generalization of our method to variable shape parameters (i.e., $\veps_i :=\veps(\bx_i)$) is straightforward.}
In this work, we only consider infinitely differentiable global basis function $\varphi^\veps$ as discussed in (\ref{GA})--(\ref{BRBF}). 
Note that the derivation and framework of our method remain the same for all basis functions listed in (\ref{GA})--(\ref{BRBF}). 

For easy explanation, we will separate our discussion of $\ap(\bx) < 2$ and $\ap(\bx) = 2$.  
For $\ap(\bx) < 2$, starting with its  integral form in (\ref{integralFL}), we first rewrite the variable-order fractional Laplacian into a summation of two integrals over $\Og$ and $\Og^c$, respectively. 
Then substituting the ansatz \eqref{Sol1D} into it and taking the extended Dirichlet boundary conditions on $\Og^c$  into account, we obtain 
\bea\label{Dlaplace1}
(-\Dt)^{{\ap(\bx)}/{2}}_hu(\bx) &=& C_{d, \ap(\bx)} \bigg({\rm P.V.}\int_\Og\fl{\widehat{u}(\bx)- \widehat{u}(\by)}{|\bx - \by|^{d+\ap(\bx)}} d\by + \int_{\Og^c}\fl{\widehat{u}(\bx) - g(\by)}{|\bx - \by|^{d+\ap(\bx)}} d\by \bigg) \nn\\
&=& (-\Dt)^{{\ap(\bx)}/{2}} \widehat{u}(\bx) +  C_{d, \ap(\bx)}  \int_{\Og^c}\fl{\widehat{u}(\by) - g(\by)}{|\bx - \by|^{d+\ap(\bx)}} d\by, 
\eea
for $0 < \ap(\bx) < 2$, where $(-\Dt)^{{\ap(\bx)}/{2}}_h$ represents  numerical approximation of $(-\Dt)^{{\ap(\bx)}/{2}}$. 
It is obvious that the integral term  in (\ref{Dlaplace1}) is caused by the extended Dirichlet boundary conditions on $\Og^c$.  
While  $\ap(\bx) = 2$, we immediately get from \eqref{Sol1D}: 
\bea\label{Dlaplace2}
(-\Dt)_hu(\bx) = -\Dt\widehat{u}(\bx), \quad \ \mbox{for} \ \ \ap(\bx) = 2. 
\eea
Combining the approximation in \eqref{Dlaplace1} and \eqref{Dlaplace2} yields a unified scheme of the variable-order Laplacian $(-\Dt)^{{\ap(\bx)}/{2}}$ for $0 <\ap(\bx) \le 2$, i.e.,  
\bea\label{Eq-uniform}
(-\Dt)_h^{{\ap(\bx)}/{2}}u(\bx) =  (-\Dt)^{{\ap(\bx)}/{2}}\widehat{u}(\bx) + C_{d, \ap(\bx)}\int_{\Og^c}\fl{\widehat{u}(\by) - g(\by)}{|\bx - \by|^{d+\ap(\bx)}} d\by,\quad 
\eea
where $C_{d, \ap(\bx)}$ is defined as in (\ref{integralFL}) { for $0 < \ap(\bx) \le 2$}.  
Note that when $\ap(\bx) = 2$,  $C_{d, 2} = 0$ and thus the integral term in \eqref{Eq-uniform} vanishes. 
The scheme (\ref{Eq-uniform}) also holds for constant-order Laplacians. 

The integral term in (\ref{Eq-uniform}) is free of singularity and  can be accurately approximated using numerical quadrature rules, while the first term in (\ref{Eq-uniform}) can be analytically expressed. 
From (\ref{Sol1D}), it is straightforward to get 
\bea\label{uhat}
(-\Dt)^{\ap(\bx)/2} \widehat{u}(\bx) = \sum_{i = 1}^{\bar{N}} \lambda_i\, \underbrace{\big[(-\Dt)^{\ap(\bx)/2}\varphi^\veps(|\bx - \bx_i|)\big]}_{{\Psi}_i^{\ap(\bx)}(\bx)}, \qquad
\eea
where for notational simplicity we denote ${\Psi}_i^{\ap(\bx)}(\bx)$ as the variable-order Laplacian of radial basis function centered at $\bx_i$. 
For example, we can easily obtain
\beas
{\Psi}_i^{\ap(\bx)}(\bx) = c_\ap\left\{\begin{array}{l}
\displaystyle \,_1F_1\Big(\fl{d+\ap(\bx)}{2}; \,\fl{d}{2};\, -\veps^2|\bx-\bx_i|^2\Big), \hspace{1cm}\\
\hspace{4.65cm} \mbox{for Gaussian RBF $\varphi^\veps$ in (\ref{GA})}, \\
\displaystyle \fl{\Gamma\big(\bt+\fl{\ap(\bx)}{2}\big)}{\Gamma(\bt)} \,_2F_1\Big(\fl{d+\ap(\bx)}{2}, \bt+\fl{\ap(\bx)}{2}; \, \fl{d}{2}; \, -\veps^2|\bx-\bx_i|^2\Big), \\
\hspace{4.65cm} \mbox{for gIMQ RBF $\varphi^\veps$  in (\ref{GIMQ})}, \\
\end{array}\right.
\eeas
with $c_\ap =|\veps|^{\ap(\bx)}  {2^{\ap(\bx)}\Gamma\big((d+\ap(\bx))/2\big)}/{\Gamma(d/2)}$ by Lemma  \ref{lemma3}. 
It is clear  that  our method avoids numerically approximating the hypersingular integral in the variable-order fractional Laplacian (\ref{integralFL}) by utilizing the analytical formulation of ${\Psi}_i^{\ap(\bx)}(\bx)$. 
Using quadrature rules to approximate hypersingular integrals could significantly increase the computational cost especially in high dimensions and make the implementation of RBF-based methods more complicated as special treatments are required around singularities \cite{Pang2015, Piret2013,  Rosenfeld2019}. 
Our method is free of these issues thanks to the results in Lemma \ref{lemma3}.

Now, we present the fully discretized scheme for (\ref{diffusion}). 
To this end, we choose test points $\bx_k \in \bar{\Og}$ for $1 \le k \le \bar{N}$.
Note that test points $\bx_k$ can be chosen independently from RBF center points $\bx_i$.
For test points in domain $\Og$, i.e., $\bx_k \in \Og$,  we substitute (\ref{Eq-uniform})--(\ref{uhat}) into the governing equation of \eqref{diffusion} { and obtain} the discretization scheme as: 
\bea\label{scheme1}
&& \sum_{i = 1}^{\bar{N}} \lambda_i \bigg[{\Psi}_i^{\ap(\bx_k)}(\bx_k)+C_{d, \ap(\bx_k)}\int_{\Og^c}\fl{\varphi^\veps(|\by-\bx_i|)}{|\bx_k - \by|^{d+\ap(\bx_k)}} d\by\bigg] \nn \\
&&\hspace{1cm}= f\big(\bx_k\big) + C_{d, \ap(\bx_k)}\int_{\Og^c}\fl{g(\by)}{|\bx_k-\by|^{d+\ap(\bx_k)}}d\by, \qquad \mbox{for} \ \ \bx_k \in\Og. \qquad 
\eea 
While for test points along boundary $\p\Og$,  i.e., $\bx_k \in \p\Og$,  we can directly apply  the ansatz (\ref{Sol1D}) to the boundary conditions in (\ref{diffusion})  and obtain
\bea\label{scheme2}
\qquad \sum_{i = 1}^{\bar{N}}\lambda_i \varphi^{\veps}(|\bx_k - \bx_i|) = g(\bx_k), \quad \ \ \mbox{for} \ \ \bx_k \in\p\Og.
\eea
That is, the boundary conditions are discretized only for points on $\p\Og$ (instead of over $\Og^c$) even if the fractional Laplacian is considered. 
In fact, if $\ap(\bx) < 2$ the extended boundary conditions on ${\mathbb R}^d\backslash\bar{\Og}$ have been taken into account via the integrals over $\Og^c$.  

Note that the fractional Laplacian is a nonlocal operator, and its nonlocality always leads to a full linear  system even if local  methods (e.g., finite difference/element methods in \cite{Duo2018, Duo-FDM2019, Acosta2017, Hao2021} for constant-order fractional Laplacian) are used.  
The variable-order Laplacian retains the same nonlocality as the fractional Laplacian. 
Besides,  the heterogeneity of variable-order Laplacian also introduces formidable challenges in numerical studies, as the resultant stiffness matrix generally does not have any symmetry. 
Currently,   nonlocality and heterogeneity remain as two main challenges in the study of variable-order fractional derivatives. 
Our methods can achieve spectral accuracy. 
They can greatly save the storage and computational cost in simulating problems with the variable-order Laplacian, especially in high-dimensional cases. 
This suggests that global methods might be more beneficial for solving nonlocal or fractional problems. 
Moreover, our methods, integrating the advantages of pseudo-differential definition (\ref{pseudo}) (i.e., compatible between $\ap = 2$ and $\ap < 2$ ) and pointwise hypersingular integral definition (\ref{integralFL}) (i.e., easy for non-periodic boundary conditions), can effectively approximate the variable-order Laplacian $(-\Dt)^{\ap(\bx)/2}$ for $0 < \ap(\bx) \le 2$.

\bigskip
In Section \ref{section4}, we will test numerical accuracy of our method, while fractional PDEs arising in various applications will be numerically studied in Section \ref{section5}. 
Unless otherwise stated, all numerical results reported here are computed by generalized inverse multiquadric function based method  with  $\bt = -(d+1)/2$ in (\ref{GIMQ}).  
Numerical studies show that the choice of $\bt$ plays an insignificant role in the performance of our method. 
More discussion on generalized inverse multiquadric basis functions can be found in \cite{Kansa90I, Kansa90II, Wu0020}. 
{ Our studies show that Gaussian RBFs can achieve similar accuracy as generalized inverse multiquadric RBFs, but Gaussian RBFs tend to require larger shape parameters. } 
In our simulations, we choose the test points $\bx_k$ from the same set of RBF center points $\bx_i$, and constant shape parameter $\veps$ is used.

\section{Approximation of variable-order Laplacians}
\label{section4}
\setcounter{equation}{0}

In this section, we test the performance of our meshfree method in approximating the variable-order Laplacian $(-\Dt)^{\ap(\bx)/2}$ for different $\ap(\bx)$. 
Here, we will focus on  the one-dimensional cases. 
Choose domain $\Og = (-1, 1)$ and consider various heterogeneous exponent $\ap(x)$, i.e.,
\bea\label{alpha-fun}
\left\{\begin{array}{l}
\ap_1(x) = 1+x,  \\
\ap_2(x) = 1-|x|,  \\
\ap_3(x) = 0.7e^{-x}, \\
\ap_4(x) = 1 + {\rm tanh}(4x+2), \\
\ap_5(x) = \cos(x), 
\end{array}\right. \qquad \mbox{for} \ \ x\in\Og.
\eea 
In the following examples, we will approximate the operator $(-\Dt)^{\ap(x)/2}$ on domain $\Og$ with different Dirichlet boundary conditions. 
Numerical errors are computed as the root mean square (RMS) error, i.e.
\bea\label{rms1}
\|e_\Dt\|_{\rm rms} = \bigg(\fl{1}{K}\sum_{l = 1}^K \big|(-\Dt)^{{\ap(x_l)}/2}u(x_l) - (-\Dt)_h^{{\ap(x_l)}/2}u(x_l)\big|^2\bigg)^{1/2},
\eea
 where $(-\Dt)_h^{\ap(x)/2}$ represents the numerical approximation of  variable-order Laplacian, and $K \gg \bar{N}$ denotes the total number of interpolation points on $\bar{\Og}$.  
Here, we choose $K$ to be large enough such that its value does not affect the RMS error in (\ref{rms1}). 
Moreover,  the choice of interpolation points $x_l$ is independent of RBF center points $x_i$ and test points $x_k$.

\bigskip 
\noindent{\bf Example 1 (Nonhomogeneous boundary conditions). }  We study and compare the performance of our method in approximating the variable-order and constant-order Laplacians.  
Choose function { $u(x) = \sqrt{2} \sin(|x|)/\big(\sqrt{\pi}|x|\big)$} for $x \in {\mathbb R}$, where the value of $u$ at $x = 0$ is defined in the limit sense (i.e., $u(0) = \sqrt{2/\pi}$). 
Here, we apply our method to numerically estimate $(-\Dt)^{\ap(x)/2} u(x)$  for $x \in \Og$, { while its exact solution is given by}
\[
{ (-\Dt)^{\ap(x)/2} u(x) = \fl{\sqrt{2}}{(\ap+1)\sqrt{\pi}}\,_1F_2\Big(\fl{1+\ap(x)}{2};\, \fl{3+\ap(x)}{2}, \fl{1}{2}; \, -\fl{x^2}{4}\Big).}
\]
Since function $u$ is defined on ${\mathbb R}$,  this problem can be viewed as approximating the Laplace operator on domain $\Og$ with extended nonhomogeneous Dirichlet boundary condition { $g(x) = \sqrt{2} \sin(|x|)/\big(\sqrt{\pi}|x|\big)$} for $x \in \Og^c$.   

Table \ref{Table5-1} presents numerical errors $\|e_\Dt\|_{\rm rms}$ for different exponent $\ap(x)$ and number of points $\bar{N}$, where shape parameter $\veps = 1$, and RBF test/center points are chosen to be uniformly distributed on $[-1, 1]$. 
It shows that as the number of points $\bar{N}$ increases,  numerical errors  decrease with a spectral rate. 
\begin{table}[htb!]
\centering
\begin{tabular}{|c||c|c|c|c|}
\hline
 \rule{0pt}{12pt} $\ap(x)$  & $\bar{N} = 5$ & $\bar{N} = 9$ & $\bar{N} = 17$ & $\bar{N} = 33$ \\ 
\hline 
 \hline
$\ap_1$ & 6.7660e-2 &2.5903e-2 & 1.4718e-3 & 1.3326e-6 \\
\hline  
$\ap_2$ &3.9105e-3 &2.9816e-4 & 2.2903e-6 & 6.686e-10 \\
 \hline 
$\ap_3$ & 3.9108e-2  & 1.4876e-2 &  8.5310e-4  & 7.6271e-7 \\
 \hline 
$\ap_4$ &9.2343e-2  & 3.1023e-2  & 1.6061e-3  & 1.3948e-6  \\
    \hline 
$\ap_5$ &7.9705e-3   &9.1226e-4  & 1.1420e-5 &  2.3138e-9  \\
    \hline
$\ap\equiv0.4$ &4.1548e-3&   4.5222e-4  & 6.1030e-6&   1.1654e-9\\
    \hline 
$\ap \equiv 1.0$ &1.3176e-2&   2.1256e-3&  4.4558e-5 &  1.3697e-8 \\
    \hline 
$\ap\equiv2.0$ &1.2938e-1 &  4.3848e-2 &  2.2715e-3 &  1.9727e-6  \\
\hline
\end{tabular}
\caption{Numerical errors  $\|e_\Dt\|_{\rm rms}$ in approximating function $(-\Dt)^{\ap(x)/2}u
(x)$ on $(-1, 1)$, where { $u(x) = \sqrt{2}\sin(|x|)/\big(\sqrt{\pi}|x|\big)$} for $x \in {\mathbb R}$, shape parameter $\veps = 1$, and $\ap_m(x)$ (for $1 \le m \le 5$) are defined in (\ref{alpha-fun}).}\label{Table5-1}
\end{table}
For the same $\bar{N}$,  numerical errors tend to be smaller if $\sup\ap(x) \le 1$ (e.g., $\ap_2$ and $\ap_5$). 
For each $\ap(x)$,  maximum numerical errors are found around  boundary points $x = \pm 1$,  consistent with our observations in RBF interpolation to function $u(x)$. 
The function $(-\Dt)^{\ap(x)/2}_hu(x)$ is computed by first obtaining coefficients $\lambda_i$ by assuming $\widehat{u}(x_k) = u(x_k)$ at all test points $x_k$, and then substituting $\lambda_i$ into (\ref{Eq-uniform}) together with (\ref{Sol1D}) and (\ref{uhat}).  
Numerical errors introduced in obtaining $\lambda_i$ are independent of exponent $\ap(x)$. 
One could further improve the accuracy by reducing errors in interpolating function $u$  \cite{Fornberg2002}. 
We find that if function $u$ is smooth enough, our method can achieve  very small numerical errors even with a small number of points, and both Gaussian or generalized inverse multiquadric RBFs yield similar numerical errors. 

In addition, Table \ref{Table5-1-1} presents numerical errors of the finite difference method to compare with our results in Table \ref{Table5-1}.  
Here, we generalize the finite difference method in \cite{Duo2018} to approximate the variable-order operator, which was originally proposed for the constant-order fractional Laplacian $(-\Dt)^\fl{\ap}{2}$ with $\ap \in (0, 2)$.
\begin{table}[htb!]
\centering
\begin{tabular}{|c||c|c|c|c|c|}
\hline\rule{0pt}{12pt}
$\ap(x)$ &    $\bar{N} = 17$ & $\bar{N} = 33$ & $\bar{N} = 65$ & $\bar{N} = 129$ \\
\hline 
$\ap_1$ & 6.6056e-5 & 1.6238e-5 & 3.9539e-6 & 9.5946e-7 \\
\hline
$\ap_2$ & 1.3482e-5 & 3.0416e-6 & 7.2571e-7  & 1.7767e-7 \\
\hline
$\ap_3$ & 3.6890e-5 & 9.1768e-6 & 2.2399e-6 & 5.4221e-7 \\
\hline
$\ap_4$ & 1.9885e-4 & 4.9608e-5 & 1.2341e-5 & 3.0670e-6 \\
\hline
$\ap_5$ & 2.2105e-5 & 5.0398e-6 & 1.2016e-6 & 2.9341e-7 \\
\hline
$\ap \equiv 0.4$ & 3.3963e-6 & 7.7919e-7 & 1.8992e-7 & 4.7775e-8\\
\hline
$\ap \equiv 1.0$ & 2.7650e-5 & 6.3336e-6 & 1.5096e-6 &  3.6811e-7 \\
 \hline
 \end{tabular}
\caption{The $l_2$-norm errors of the finite difference method \cite{Duo2018} in approximating function $(-\Dt)^{\ap(x)/2}u(x)$ on $(-1, 1)$, where { $u(x) = \sqrt{2}\sin(|x|)/\big(\sqrt{\pi}|x|\big)$} for $x \in {\mathbb R}$, mesh size $h = 2/(\bar{N}-1)$, and $\ap_m(x)$ (for $1 \le m \le 5$) are defined in (\ref{alpha-fun}).}\label{Table5-1-1}
\end{table}
It is clear that to achieve the same accuracy, our RBF-based method requires much fewer points. 
Moreover, our method allows $\ap(\bx) = 2$, while finite difference method requires $0 < \ap(\bx) < 2$ since it is developed based on the integral definition in (\ref{integralFL}).

\bigskip 
\noindent{\bf Example 2 (Homogeneous boundary conditions). }  
In this example, we consider compactly supported  functions  on $\Og$. 
This problem can be viewed as approximating the Laplacian with extended homogeneous Dirichlet boundary conditions. 
Choose function $u(x) = \big[(1-x^2)^s\big]_+$ for constant $s > 0$. 
Different from $u \in C^\infty({\mathbb R})$ in Example 1, function $u$ in this case has finite smoothness at points $x = \pm 1$, depending on the value of $s$. 

Table \ref{Table5-2} shows numerical errors in approximating $(-\Dt)^{\ap(x)/2}u$ with $u(x) = (1-x^2)_+$ for different $\ap(x)$ and  $\bar{N}$. 
{ The exact solution of $(-\Dt)^{\ap(x)/2} u(x)$ is given by}
\[
{ (-\Dt)^{\ap(x)/2} u(x) = \fl{2^{\ap(x)}\Gamma\big(\fl{1+\ap(x)}{2}\big)}{\sqrt{\pi}\Gamma\big(2-\fl{\ap(x)}{2}\big)}\,_2F_1\Big(\fl{1+\ap(x)}{2},  -1+\fl{\ap(x)}{2};\, \fl{1}{2}; \, x^2\Big),}
\]
{ for  $-1 < x < 1$. }
In our simulations, we choose RBF center/test points  uniformly distributed on $[-1, 1]$.  
From Tables \ref{Table5-1} and \ref{Table5-2}, we find that both smoothness of  $u(x)$ on $\bar{\Og}$ and value of  $\ap(x)$ have impacts on numerical accuracy in approximating the Laplacian operators.  
Numerical errors are generally larger if $\sup\ap(x) > 1$, a similar observation as in Table \ref{Table5-1}. 
\begin{table}[htb!]
\centering
\begin{tabular}{|c||c|c|c|c|}
\hline
 \rule{0pt}{12pt} $\ap(x)$ & $\bar{N} = 5$ & $\bar{N} = 9$ & $\bar{N} = 17$ & $\bar{N} = 33$ \\
\hline  
\multirow{1}{*}{$\ap_1$} & 1.0027 & 4.4181e-1 & 3.9168e-2 & 9.1883e-5 \\
\hline     
\multirow{1}{*}{{$\ap_2$}}  & 5.1815e-2 & 4.6255e-3 & 5.8942e-5 & 2.9064e-8 \\
\hline    
\multirow{1}{*}{$\ap_3$} & 5.7334e-1 & 2.5332e-1 & 2.2695e-2 & 5.2574e-5 \\
\hline 
\multirow{1}{*}{$\ap_4$}  & 1.3307 & 5.2260e-1 & 4.2654e-2 & 9.6097e-5 \\
 \hline 
\multirow{1}{*}{{$\ap_5$}}  & {1.0689e-1}&{1.4340e-2} &{2.9488e-4} &{1.5848e-7} \\
\hline 
$\ap \equiv 0.4$  & 5.6460e-2  &  7.2266e-3 &   1.5830e-4 &   8.1212e-8 \\
 \hline 
$\ap \equiv 1.0$  & 1.8008e-1  & 3.4163e-2  & 1.1604e-3   & 1.0146e-6\\
\hline  
$\ap \equiv 1.5$  & 5.3719e-1   & 1.3938e-1  & 7.1335e-3  & 1.0443e-5\\
\hline
\end{tabular}
\caption{Numerical errors  $\|e_\Dt\|_{\rm rms}$ in approximating function $(-\Dt)^{\ap(x)/2}u
(x)$ on $(-1, 1)$, where $u = (1-x^2)_+$,  shape parameter $\veps = 1$, and $\ap_m(x)$ (for $1 \le m \le 5$) are defined in (\ref{alpha-fun}).}\label{Table5-2}
\end{table}
Compared to  nonhomogeneous boundary conditions, computations in this example take shorter time as the integral of $g(x)$ over $\Og^c$ is always zero. 
But, the boundary conditions do not affect the accuracy of our method. 
%
\section{Solutions of variable-order fractional PDEs}
\label{section5}
\setcounter{equation}{0}

In this section, we apply our method to study the solutions of fractional PDEs with variable-order Laplacian. 
The variable-order Laplacian makes the study of heterogenous media much easier, but it also introduces considerable challenges in numerical simulations. 
To the best of our knowledge, no numerical method has been reported for this operator. 
The lack of numerical methods greatly hinders the application of variable-order Laplacian and study of heterogenous media. 
We provide the first numerical methods for the variable-order Laplacian $(-\Dt)^{\ap(\bx)/2}$. 
Our methods can not only solve  problems with the fractional Laplacian for $0 < \ap(\bx) < 2$, but also allow to include exponent $\ap(\bx) = 2$,  providing a unified scheme for classical and fractional Laplacians. 
This  important property enables us to easily study the coexistence and transition of anomalous and normal diffusion in many complex systems. 

\subsection{Poisson problems} 
\label{section5-1}

So far, the constant-order fractional Poisson problems with extended homogeneous Dirichlet boundary conditions (i.e., $g(\bx) \equiv 0$) have been well studied both analytically and numerically.  
However, the understanding of variable-order fractional Poisson problems  still falls very behind. 
To the best of our knowledge, no numerical results have been reported even for this simplest variable-order problem. 
Here,  we consider the one-dimensional (i.e., $d = 1$) Poisson problem (\ref{diffusion}) on domain $\Og = (-1, 1)$. 
We will first study the accuracy of our method and then  explore the effect of exponent $\ap(x)$ on solution of Poisson equation.  

Table \ref{Table5-3} presents numerical errors in solution $u$ for different exponent $\ap(x)$ as listed in (\ref{alpha-fun}), where  $f$ is chosen such that the exact solution of (\ref{diffusion}) is $u(x) = \big[(1-x^2)^3\big]_+$ for $x \in {\mathbb R}$. 
\begin{table}[htb!]
\centering
\begin{tabular}{|c||c|c|c|c|c|}
\hline
\rule{0pt}{12pt} $\ap(x)$ &    $\bar{N} = 5$ & $\bar{N} = 9$ & $\bar{N} = 17$ & $\bar{N} = 33$ & $\bar{N} = 65$ \\
\hline 
    $\ap_1$ & 6.4786e-2 &  1.1510e-3  & 6.0664e-4 &  6.4472e-5  & 2.9989e-7\\ 
    \cline{2-6}
    \hline 
    $\ap_2$ &   6.4359e-2  & 9.4464e-4  & 1.6883e-4 &  8.4450e-6   &1.4821e-8\\ 
    \cline{2-6}
    \hline  
    $\ap_3$&  5.2901e-2  & 1.0656e-3 &  4.4010e-4  & 4.3140e-5  & 1.6498e-7\\ 
    \cline{2-6}
    \hline 
    $\ap_4$ &  1.8067e-1 &  2.7515e-3 &  9.3588e-4  & 1.0025e-4  & 2.4128e-7\\ 
    \cline{2-6}
    \hline 
    $\ap_5$ & 6.6009e-2 &  1.0023e-3  & 2.4019e-4 &  1.3242e-5  & 7.1019e-8\\ 
    \cline{2-6}
    \hline
    $\ap \equiv 0.4$ &  4.6934e-2 &  9.2165e-4 &  1.8985e-4  & 1.0400e-5  & 2.1501e-8\\ 
    \cline{2-6}
    \hline 
    $\ap \equiv1.0$ &  6.6192e-2 &  1.0330e-3 &  3.3376e-4  & 2.1650e-5  & 6.9136e-8\\ 
    \cline{2-6}
    \hline 
    $\ap \equiv2.0$  &2.0163e-1  & 3.5564e-3  & 1.6341e-3  & 1.7738e-4 &  6.9761e-7\\ 
    \cline{2-6}
    \hline
    \end{tabular}
\caption{The RMS errors $\|e_u\|_{\rm rms}$ of our method in solving the Poisson problem \eqref{diffusion} with exact solution $u(x) = \big[(1-x^2)^3\big]_+$, where the shape parameter $\veps = 2$.}\label{Table5-3}
\end{table}
In our simulations,  we set the shape parameter $\veps = 2$ and choose RBF center/test points uniformly distributed on $\bar{\Og}$. 
The root mean square errors in solution $u$ are calculated as
\bea\label{rms2}
\|e_u\|_{\rm rms} = \bigg(\fl{1}{K}\sum_{l = 1}^K \big|u(x_l) - \widehat{u}(x_l)\big|^2\bigg)^{1/2},
\eea
where $u$ and $\widehat{u}$ represents the exact and numerical solutions, respectively, and $K \gg \bar{N}$ denotes the total number of interpolation points on $\bar{\Og}$. 
It shows that numerical errors decrease quickly with the number of points increasing. 
Moreover, the implementations of our method for variable-order and constant-order Laplacians are  the same. 
Generally, the spatial dependence of $\ap(x)$ destroys the symmetry of stiffness matrix, and thus one has to save the entire dense matrix in simulations. 
This considerably increases the storage and computational costs, especially if low-order numerical methods are used. 
However, our RBF-based methods can achieve higher accuracy with much less points. 
\begin{table}[htb!]
\centering
\begin{tabular}{|c||c|c|c|c|c|c|}
\hline\rule{0pt}{12pt}
$\ap(x)$ &    $\bar{N} = 17$ & $\bar{N} = 33$ & $\bar{N} = 65$ & $\bar{N} = 129$ & $\bar{N} = 257$ & $\bar{N} = 513$ \\
\hline 
$\ap_1$ & 3.6456e-3 &  8.9956e-4 &  2.1469e-4  & 5.0700e-5  & 1.1951e-5  & 2.8220e-6 \\
$\ap_2$ &9.0142e-4  & 2.0260e-4  & 4.8621e-5  & 1.1942e-5  & 2.9613e-6 &  7.3747e-7\\
$\ap_3$ &2.4340e-3  & 5.6596e-4  & 1.2901e-4  & 2.9566e-5  & 6.8602e-6  & 1.6129e-6\\
$\ap_4$ &8.3653e-3  & 2.0015e-3  & 4.8355e-4   &1.1732e-4  & 2.8542e-5  & 6.9598e-6\\
$\ap_5$ &9.7335e-4  & 2.2776e-4  & 5.5704e-5   &1.3714e-5  & 3.3890e-6  & 8.4061e-7\\
{ $\ap \equiv 0.4$} & { 7.7547e-4} &  { 2.2030e-4}  &  { 5.7436e-5} &  { 1.4560e-5}  &  { 3.6569e-6}   & { 9.1562e-7}\\
{ $\ap \equiv 1.0$} & { 6.9118e-4}  &{ 9.4901e-5} &{  1.3531e-5} &{   2.3628e-6} &{  5.2722e-7} &{  1.3072e-7} \\
 \hline
\end{tabular}
\caption{The $l_2$-norm errors of the finite difference method \cite{Duo2018} in solving the Poisson problem \eqref{diffusion} with exact solution $u(x) = \big[(1-x^2)^3\big]_+$, where mesh size $h = 2/(\bar{N}-1)$.}\label{Table5-3-1}
\end{table}
To see this, we present in Table \ref{Table5-3-1} the numerical errors of the finite difference method  \cite{Duo2018} to compare with the results in Table \ref{Table5-3}. 
It shows that the finite difference method requires more points to achieve the same accuracy. 
{ Moreover, finite difference method discretizes the integral form in (\ref{integralFL}) and thus requires $0 < \ap(x) < 2$.}
Our method is more advantageous to study problems with variable-order Laplacian { $0 < \ap(x) \le 2$}. 

Next, we study the solution of fractional Poisson problem to further explore the  heterogeneous effects of $\ap(x)$, where we choose $f(x) = 2\sin^2(\pi x)$ and $g(x) \equiv 0$ in (\ref{diffusion}).
\begin{figure}[htb!]
\centerline{
\includegraphics[height = 4.86cm, width = 6.70cm]{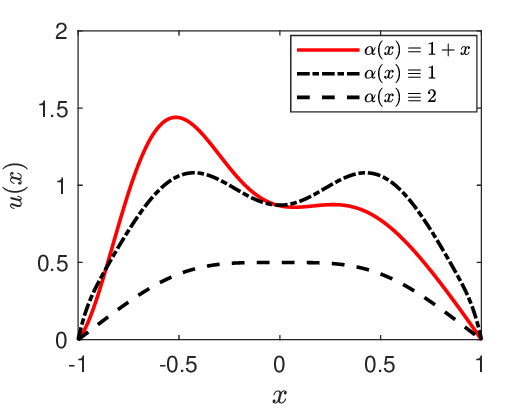}\hspace{-2mm}
\includegraphics[height = 4.86cm, width = 6.70cm]{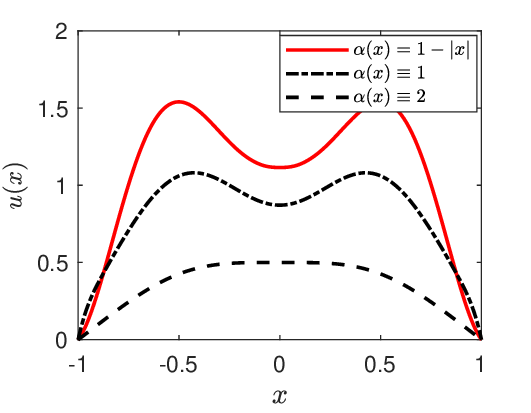}} 
\vspace{-1mm}
\centerline{
\includegraphics[height = 4.86cm, width = 6.70cm]{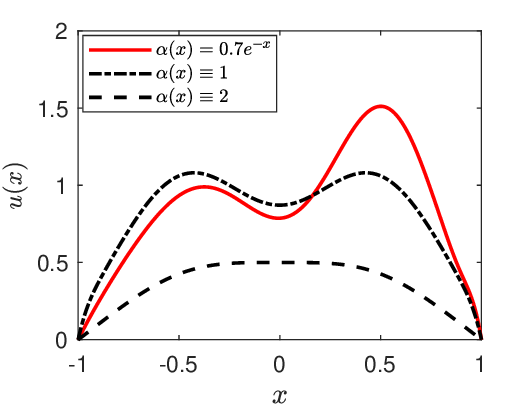}\hspace{-2mm}
\includegraphics[height = 4.86cm, width = 6.70cm]{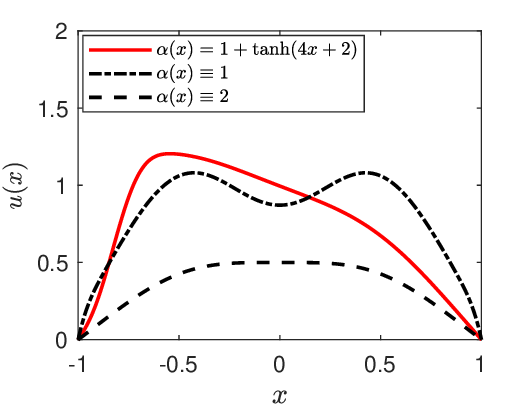}
}
\caption{Numerical solution of the one-dimensional Poisson problem in (\ref{diffusion}) with different $\ap(x)$, where $f(x) = 2\sin^2(\pi x)$ and $g(x) \equiv 0$.}\label{Figure5-3-1} 
\end{figure}
In this case, the exact solution of fractional Poisson problem is unknown. 
Figure \ref{Figure5-3-1} presents numerical solution for different $\ap(x)$. 
For ease of comparison, we also include the solution of Poisson problem with classical { negative} Laplacian (i.e., $\ap(x) \equiv 2$) or square root of the { negative}  Laplacian (i.e., $\ap(x) \equiv 1$). 
Since $f$ is symmetric about $x = 0$, the solution retains this symmetry if $\ap(x) = \ap(-x)$. 
The heterogeneity of variable-order Laplacian has strong impacts on the solutions of Poisson problem,  which makes them significantly different from solutions of constant-order Poisson problems. 
\subsection{Wave propagation in heterogeneous media}
\label{section5-2}

Wave propagation in heterogeneous media has been widely studied with fractional PDEs models, where different modeling approaches have been used to describe the heterogeneous media; see \cite{Meerschaert2008, Zhu2014, Xue2018, Chen2014b} and references therein. 
Here, we consider the one-dimensional fractional wave equation of the form \cite{Zhu2014, Xue2018, Chen2014b}: 
\bea\label{wave}
\fl{1}{c^2}\p_{tt} u(x, t) = - (-\Dt)^{\ap(x)/2} u, \qquad {\mbox{for} \  t\in (0, T]}
\eea
with constant $c  > 0$. 
The initial conditions are taken as 
\bea\label{wave-IC}
u(x, 0) = {\rm sech}\big(a(x+2)\big),\quad u_t(x, 0) = b\,{\rm sech}\big(a(x+2)\big){\rm tanh}\big(a(x+2)\big).\ 
\eea
In the classical case with $\ap(x) \equiv 2$, the wave equation (\ref{wave})--(\ref{wave-IC}) admits the exact solution of the form
\bea\label{ex}
u_{\rm ex}(x, t) = {\rm sech}\big[a(x+2)-b\,t\big], \qquad \mbox{for} \ \ x\in {\mathbb R}, \quad t \ge 0,
\eea
provided that  $c^2 a^2 = b^2$. 
However, the exact solution in fractional cases still remains unknown even for (\ref{wave}) with constant-order fractional Laplacian $(-\Dt)^\fl{\ap}{2}$. 
{ In existing studies (e.g., \cite{Zhu2014, Xue2018, Chen2014b}), the variable-order fractional Laplacian is approximated by its constant-order counterparts. 
This reduces the computational complexity in numerical studies but also alters the heterogeneity of the model.} 

\begin{figure}[htb!]
\centerline{
(a)\includegraphics[height = 4.8cm, width = 6.70cm]{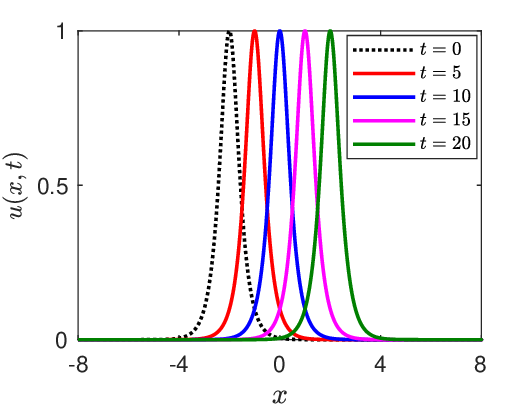}\hspace{-5mm}
(b)\includegraphics[height = 4.8cm, width = 6.70cm]{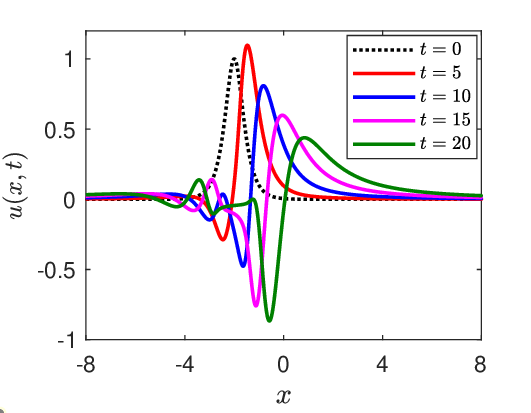}} 
\centerline{
(c)\includegraphics[height = 4.8cm, width = 6.70cm]{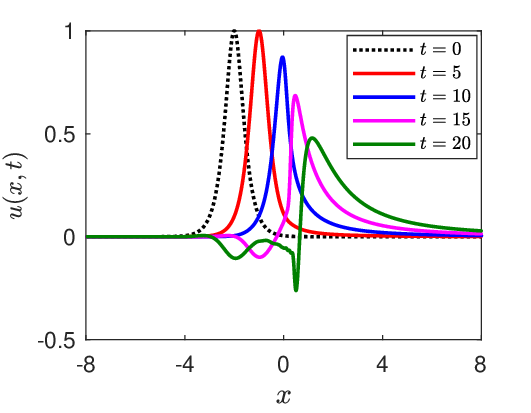}\hspace{-5mm}
(d)\includegraphics[height = 4.8cm, width = 6.70cm]{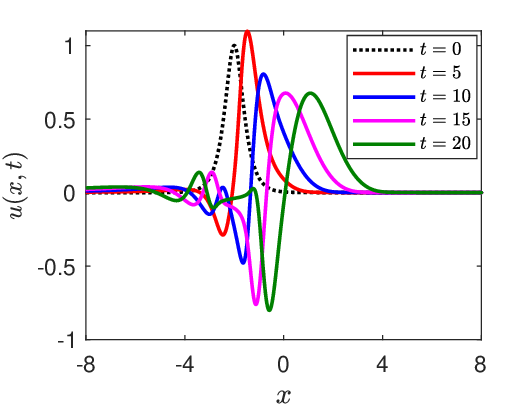}}
\caption{Time evolution of wave solutions for different $\ap(x)$, where  (a) $\ap(x) \equiv 2$; (b) $\ap(x) \equiv 1.2$; (c) $\ap(x) =  1.6-0.4{\rm tanh}(5x)$; (d) $\ap(x) =  1.6+0.4{\rm tanh}(5x)$.  
For better illustration, we only display the solution on  interval $[-8, 8]$,  much smaller than the actual computational domain. }\label{Figure6-1-1}
\end{figure}

In our simulations, {we choose $T = 20$},  and the computational domain is taken as $\Og = (-20, 20)$ with  extended homogeneous Dirichlet boundary conditions on $\Og^c$.  
It has been verified that  $\Og$ is large enough  such that the effects of domain truncation can be neglected in our simulations. 
The spatial discretization of wave equation (\ref{wave}) is realized by our meshfree method, while the time is discretized by central difference scheme with  step  size $\tau = 0.001$.  
We choose RBF center/test points uniformly distributed on $\bar{\Og}$ with $\bar{N} = 641$, and the shape parameter $\veps = 2$. 
Figure \ref{Figure6-1-1} shows the time evolution of  wave solution for different exponent $\ap(x)$, where we choose parameters $c = 0.2, a = 3$, and $b = 0.6$ in (\ref{wave})--(\ref{wave-IC}).
The soliton-like solution is initially centered at $x = -2$, and then it travels from left to right over time. 
In classical media with $\ap(x) \equiv 2$, the shape of solution $u$ remains the same for any time $t \ge 0$ (see Figure \ref{Figure6-1-1} (a)), consistent with the exact solution in (\ref{ex}).

In contrast to classical cases,  solutions of the fractional wave equation lose its original soliton-like shape, and  scattering of waves is observed even in constant-order fractional cases (see Figure \ref{Figure6-1-1} (b)).  
These phenomena are similar to solution decoherence observed in the fractional nonlinear Schr\"odinger equations \cite{Duo2016, Kirkpatrick2016}. 
In heterogeneous cases, we choose $\ap(x) = 1.6 \pm 0.4{\rm tanh}(5x)$ and study wave transition between classical and fractional media.  
Here, the transition between two different media is described by a hyperbolic tangent function  \cite{Zhu2014}.  
Figure \ref{Figure6-1-1} (c) illustrates the wave traveling from classical to fractional media. 
It shows that the solution behaves like in classical media before the wave front reaches transition region. 
Once entering the fractional media, the wave shape is distorted, and radiations of waves  due to the nonlocality of fractional Laplacian are observed. 
Figure \ref{Figure6-1-1} (d) shows that from fractional to classical media, the soliton-like solution distorts and radiates from the beginning, and  the original shape cannot be restored even after the wave enters classical media. 
Moreover,  radiations of waves are only observed in  fractional media but not in classical media. 
Future studies will be carried out to further understand wave propagation in heterogeneous media.

\subsection{Coexistence of normal and anomalous diffusion}
\label{section5-3}

Recent literature shows that normal and anomalous diffusion may coexist in many
complex systems \cite{Javanainen2013, Zhang2012, Lenzi2016}. 
The variable-order Laplacian makes it much easier to study such a coexistence by controlling the spatial-dependent exponent $\ap(\bx)$. 
In this section, we compare normal and anomalous diffusion and study the coexistence of these two  diffusion processes. 
Let domain $\Og = \Xi\backslash\big(\Xi_+ \cup \Xi_-\big) \in {\mathbb R}^2$ be an irregular channel  with  $\Xi = (-3, 3)\times(-1, 1)$ and $\Xi_\pm = [-1, 1] \times (\pm1, \pm0.5]$. 
We consider the following diffusion problem:
\bea\label{diff}
\p_t u(\bx, t) = {-\kappa} (-\Delta)^{{\ap(\bx)}/{2}} u,  \qquad\mbox{for} \ \ \bx \in \Og, \quad t > 0
\eea
with extended homogeneous Dirichlet boundary conditions, where $\kappa > 0$ denotes the diffusion coefficient. 
The initial condition is taken as 
\bea\label{diff-ic}
u(\bx, 0) = \left\{\begin{array}{ll} 
1, \quad \ & \mbox{if} \ \ \bx \in [-0.5, 0.5]^2,  \\
0, & \mbox{otherwise,}
\end{array}\right. \qquad \mbox{for} \ \ \bx \in \bar{\Og},
\eea
that is, initially the solution concentrates on a square region $[-0.5, 0.5]^2$ at the center of channel. 
See Figure \ref{Figure6-2-4} ($t = 0$) for the illustration of domain $\Og$ and initial condition $u(\bx,0)$.  
In our simulations, we use the Crank--Nicolson method for time discretization with step size $\tau = 0.001$. 
RBF center/test points are chosen as equally spaced grid points on $\bar{\Og}$ with $\bar{N} = 713$, and the shape parameter is set as $\veps = 2$. 

\begin{figure}[htb!]
\centerline{
\includegraphics[height = 2.6cm, width = 4.76cm]{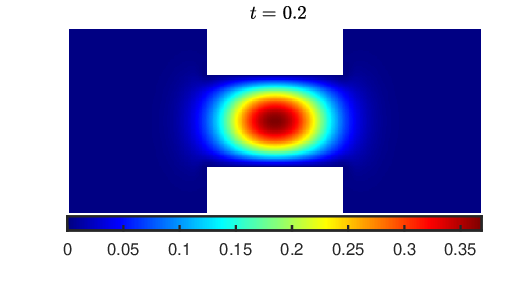}\hspace{-5mm}
\includegraphics[height = 2.6cm, width = 4.76cm]{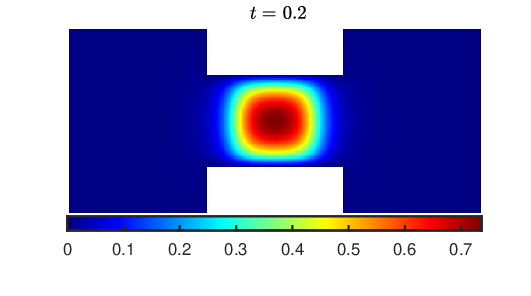}\hspace{-5mm}
\includegraphics[height = 2.6cm, width = 4.76cm]{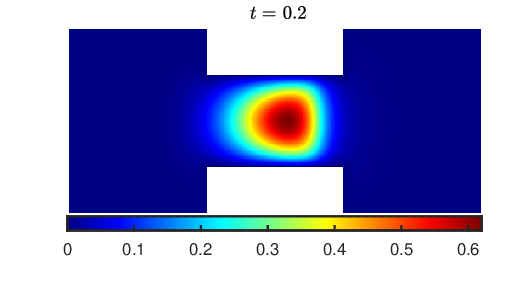}}\vspace{0mm}
\centerline{
\includegraphics[height = 2.6cm, width = 4.76cm]{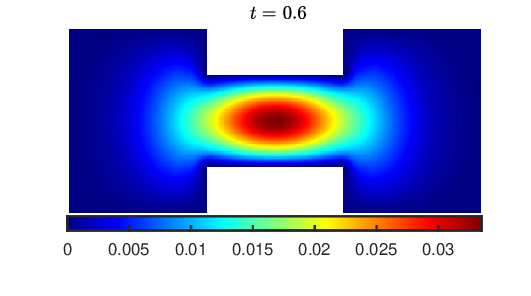}\hspace{-5mm}
\includegraphics[height = 2.6cm, width = 4.76cm]{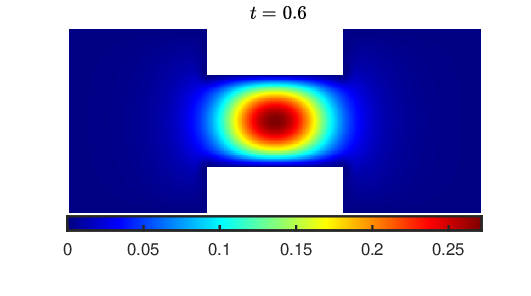}\hspace{-5mm}
\includegraphics[height = 2.6cm, width = 4.76cm]{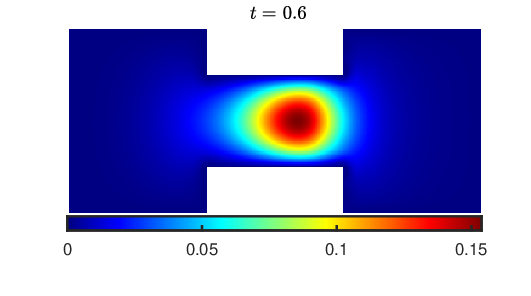}}\vspace{0mm}
\centerline{
\includegraphics[height = 2.6cm, width = 4.76cm]{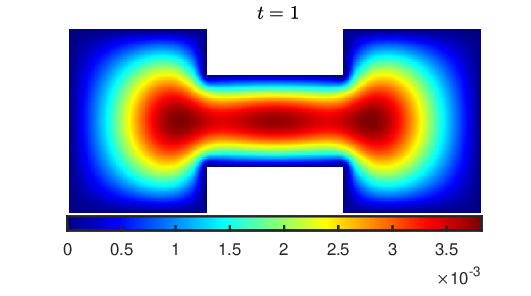}\hspace{-5mm}
\includegraphics[height = 2.6cm, width = 4.76cm]{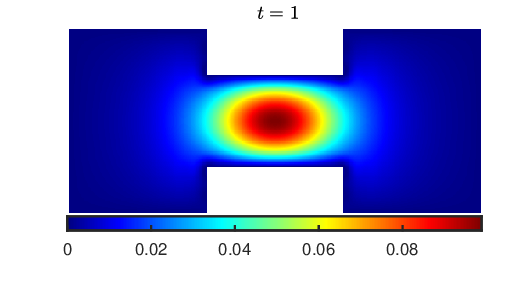}\hspace{-5mm}
\includegraphics[height = 2.6cm, width = 4.76cm]{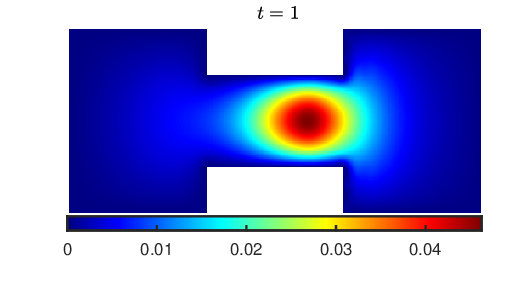}}\vspace{0mm}
\centerline{
\includegraphics[height = 2.6cm, width = 4.76cm]{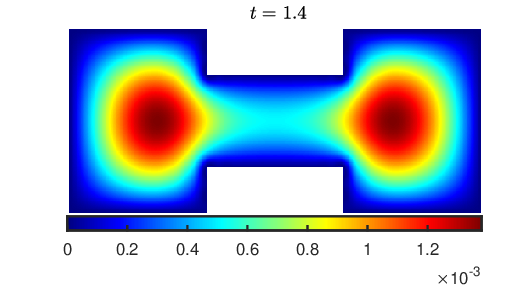}\hspace{-5mm}
\includegraphics[height = 2.6cm, width = 4.76cm]{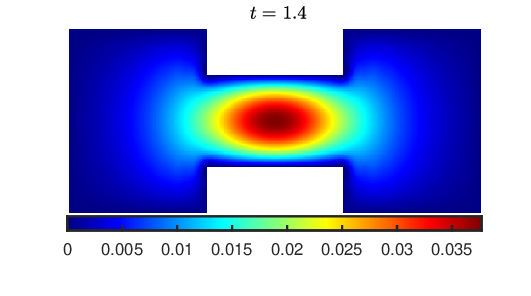}\hspace{-5mm}
\includegraphics[height = 2.6cm, width = 4.76cm]{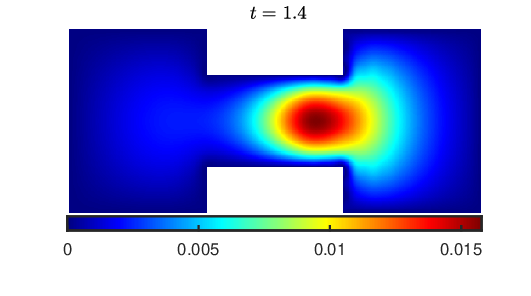}}\vspace{0mm}
\centerline{
\includegraphics[height = 2.6cm, width = 4.76cm]{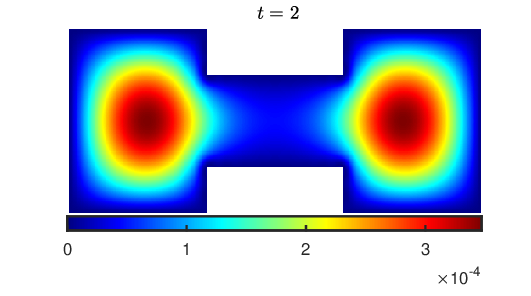}\hspace{-5mm}
\includegraphics[height = 2.6cm, width = 4.76cm]{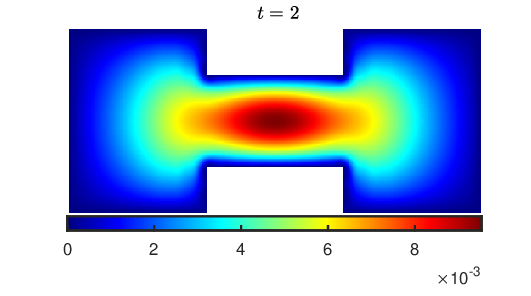}\hspace{-5mm}
\includegraphics[height = 2.6cm, width = 4.76cm]{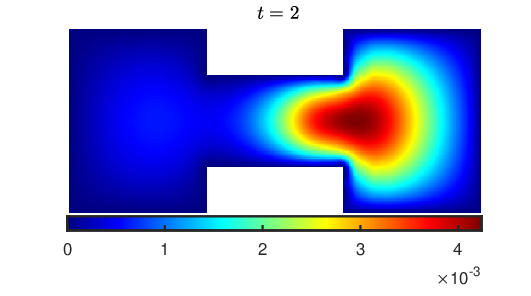}}\vspace{0mm}
\caption{Time evolution of solution $u(\bx, t)$ of (\ref{diff})--(\ref{diff-ic}), where $\ap(\bx) \equiv 2$ (left column),  $\ap(\bx) \equiv 1.4$ (middle column), and $\ap(\bx) = 2{\chi_{\{x \le -0.5\}}} + 1.4{\chi_{\{x \ge 0.5\}}} + \big(1.7-0.6x\big){\chi_{\{|x| < 0.5\}}}$ (right column). }\label{Figure6-2-2}
\end{figure}

In Figure \ref{Figure6-2-2}, we compare the solution dynamics in normal ($\ap(\bx) \equiv 2$), anomalous ($\ap(\bx) \equiv 1.4$), and their coexisting ($\ap(\bx) = 2{\chi_{\{x \le -0.5\}}} + 1.4{\chi_{\{x \ge 0.5\}}} + \big(1.7-0.6x\big){\chi_{\{|x| < 0.5\}}}$) systems, where we  choose $\kappa = 0.5$. 
During the dynamics, the solution diffuses to both sides of the channel, and at the same time it decays due to homogeneous Dirichlet boundary conditions.  
In the system with normal ($\ap(\bx) \equiv 2$) diffusion,  the solution spreads symmetrically along the $x$-axis and quickly decays through the contact with zero boundary conditions. 
For example, the solution reduces to $u \sim {O}(10^{-4})$ at time $t = 2$. 
For $\ap(\bx) \equiv 1.4$, the diffusion is still homogeneous in space, but the solution decays much slower. 
Our extensive studies show that the smaller the value of constant $\ap$, the slower the solution diffuses. 
Different from normal cases,  zero boundary conditions can affect the solution through both contact  interactions along $\p\Og$ and long-range interactions from ${\mathbb R}^2\backslash\bar{\Og}$. 

The solution evolution in heterogeneous diffusion system is very different (see Figure \ref{Figure6-2-2} right column), where normal and anomalous diffusion coexist. 
In this case, normal diffusion characterizes the region of $x \le -0.5$, while anomalous diffusion becomes dominant if $x > -0.5$. 
But, the anomalous diffusion affects the solution on whole domain  $\Og$ due to long-range interactions. 
\begin{figure}[htb!]
\centerline{
\includegraphics[height = 4.16cm, width = 5.26cm]{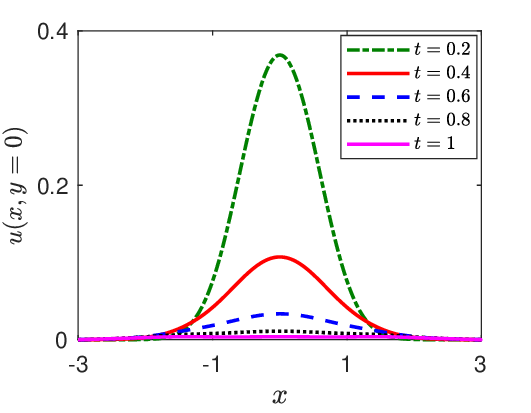}\hspace{-2mm}
\includegraphics[height = 4.16cm, width = 5.26cm]{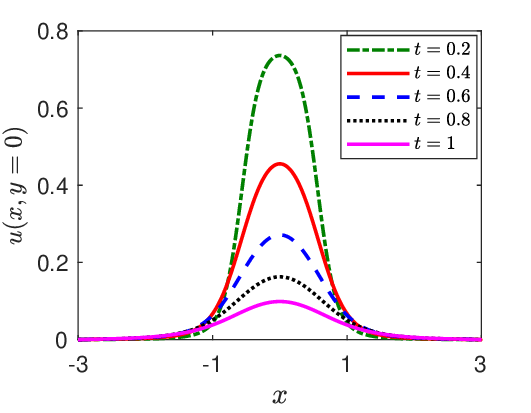}\hspace{-2mm}
\includegraphics[height = 4.16cm, width = 5.26cm]{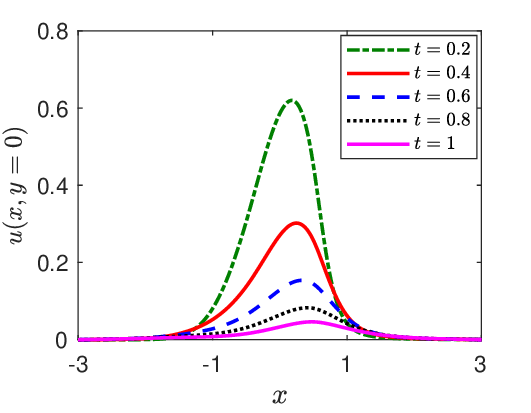}}
\caption{Time evolution of $u(\bx, t)$ at $y = 0$, where $\ap(\bx) \equiv 2$ (left column),  $\ap(\bx) \equiv 1.4$ (middle column), and $\ap(\bx) = 2{\chi_{\{x \le -0.5\}}} + 1.4{\chi_{\{x \ge 0.5\}}} + \big(1.7-0.6x\big){\chi_{\{|x| < 0.5\}}}$ (right column).}\label{Figure6-2-3}
\end{figure} 
In this case, the solution is asymmetric in $x$-direction over time $t >0$ (see illustration in Figures \ref{Figure6-2-2}--\ref{Figure6-2-3}).
Moreover, the solution of coexisting systems diffuses slower than that from $\ap(\bx) \equiv 2$, but faster than $\ap(\bx) \equiv 1.4$. 
It is clear that nonlocal interactions from variable-order Laplacian have strong impacts on the solution, and coexistence of normal and anomalous diffusion could significantly change the dynamics compared to homogeneous diffusion systems.

In Figure \ref{Figure6-2-4}, we further demonstrate the effects of heterogeneous diffusion by taking $\ap(\bx)  = (9+x)/6$. 
In this case, exponent $\ap(\bx)$  increases from $\ap(x=-3) = 1$ to $\ap(x=3) = 2$, leading to a continuous transition between anomalous diffusions. 
\begin{figure}[htb!]
\centerline{
\includegraphics[height = 2.6cm, width = 4.76cm]{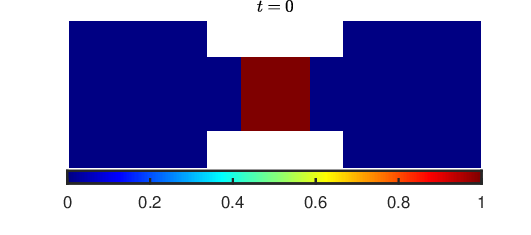}\hspace{-5mm}
\includegraphics[height = 2.6cm, width = 4.76cm]{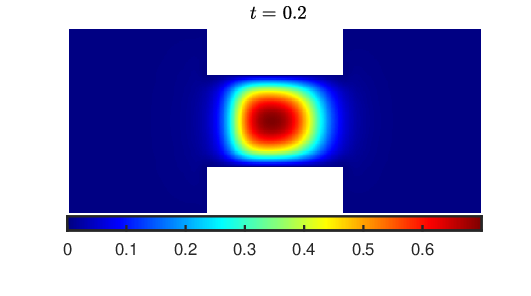}\hspace{-5mm}
\includegraphics[height = 2.6cm, width = 4.76cm]{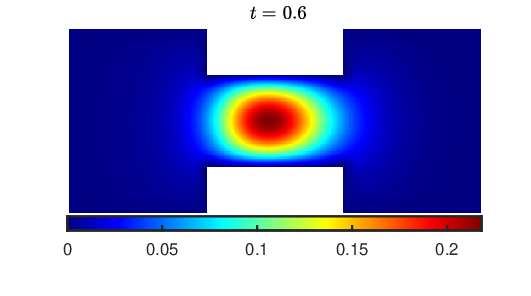}}\vspace{0mm}
\centerline{
\includegraphics[height = 2.6cm, width = 4.76cm]{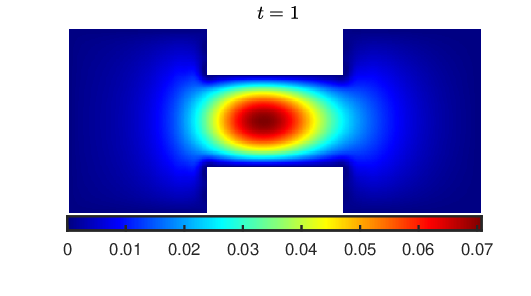}\hspace{-5mm}
\includegraphics[height = 2.6cm, width = 4.76cm]{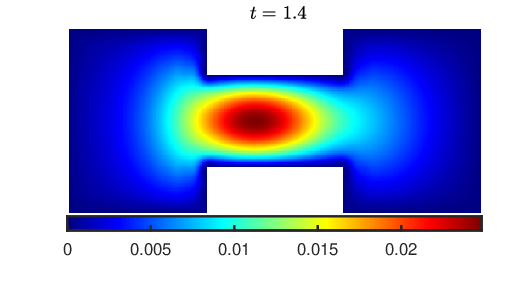}\hspace{-5mm}
\includegraphics[height = 2.6cm, width = 4.76cm]{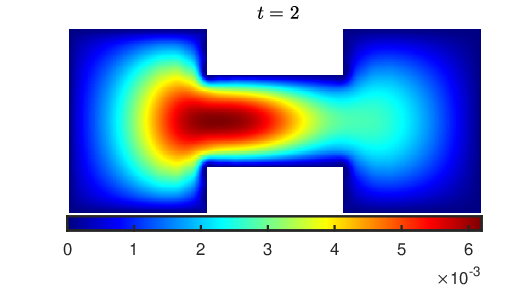}}\vspace{-2mm}
\caption{Time evolution of solution $u(\bx, t)$ of (\ref{diff})--(\ref{diff-ic}) with $\kappa = 1$ and $\ap(\bx) = (9+x)/6$ .}\label{Figure6-2-4}
\end{figure}  
Consistent with our observations in Figure \ref{Figure6-2-2}, the larger the exponent $\ap(\bx)$, the faster the solution diffuses. 
The solutions become asymmetric due to heterogeneous effects, and  exponent $\ap(\bx)$ plays an important role on solution dynamics. 
Numerical studies show that our methods are highly effective in solving problems with variable-order Laplacian and could be easily applied to study various heterogeneous problems arising in different fields \cite{Javanainen2013, Zhang2012, Lenzi2016}. 

\subsection{Allen--Cahn problems in heterogeneous fluids}
\label{section5-4}

The Allen--Cahn equation is well-known in modeling phase field problems arising in material sciences and fluid dynamics. 
Recently, its fractional analogue has been proposed to study phase transition in the presence of anomalous diffusion \cite{Song2016, Duo2019}.  
In the following, we apply our method to study  coalescence of two ``kissing" bubbles in the heterogeneous fractional Allen--Cahn equation of the form \cite{Song2016, Duo2019}: 
\begin{eqnarray}\label{AC-1}
	\p_tu(\bx, t) = -(-\Delta)^{{\ap(\bx)}/{2}}u - \frac{1}{\dt^2}\,u(u^2-1),  &\quad &\mbox{for} \ \, \bx \in \Omega, \ \ t>0,\qquad\qquad \\ \label{AC-2}
 u(\bx, t) = -1, &&\mbox{for} \ \, \bx \in \Upsilon, \ \ t\geq 0, 
\end{eqnarray}
where  $u$ is the phase field function, and $\dt > 0$ represents the diffuse interface width. 
Let the domain $\Og = (0, 1)^2$, and choose the initial condition
\begin{eqnarray}\label{u0-two_bubble}
	u(\bx, 0) = 
		1-\tanh\bigg(\frac{ |\bx -  \bx_1^0| - {0.12}}{{\dt}} \bigg) - \tanh\bigg(\frac{ |\bx -  \bx_2^0| - {0.12}}{{\dt}} \bigg).  
\end{eqnarray}
In our simulations, we choose $\delta = 0.1$, $\bx_1^0 = (0.38, 0.38)$, $\bx_2^0 = (0.62, 0.62)$. 
The time of (\ref{AC-1}) is discretized by the fourth order Runge--Kutta method with time step $\tau = 0.001$, and  RBF center/test points are chosen as equally-spaced grid points with $\bar{N} = 256$. 
The shape parameter is set as $\veps = 2$. 
		
\begin{figure}[htb!]
\centerline{
\includegraphics[height = 3.28cm, width = 3.6cm]{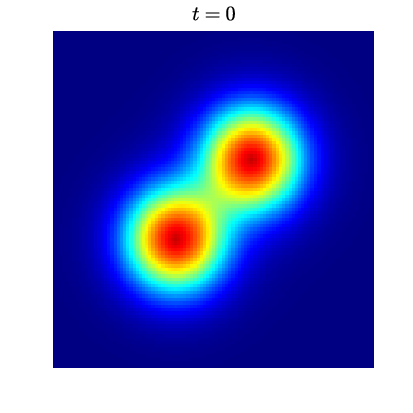}\hspace{-4mm}
\includegraphics[height = 3.28cm, width = 3.6cm]{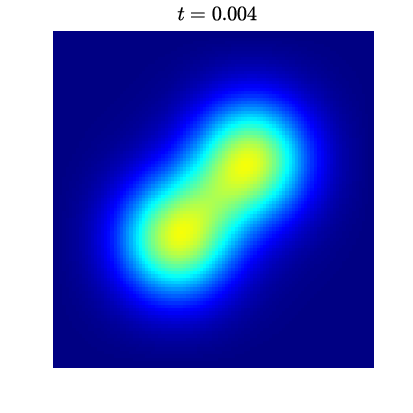}\hspace{-4mm}
\includegraphics[height = 3.28cm, width = 3.6cm]{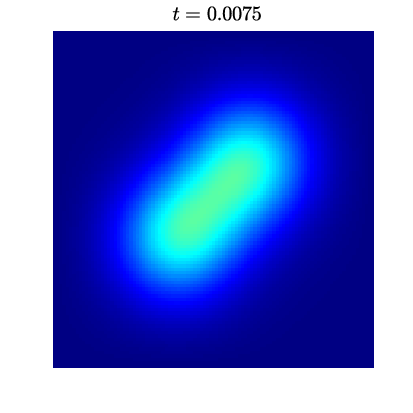}\hspace{-4mm}
\includegraphics[height = 3.28cm, width = 3.6cm]{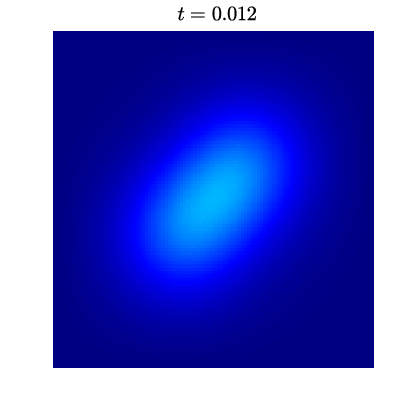}}\vspace{-2mm}
\centerline{
\includegraphics[height = 3.28cm, width = 3.6cm]{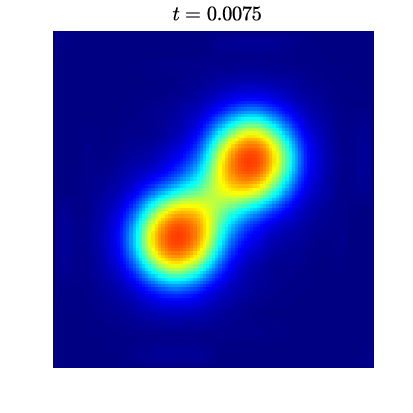}\hspace{-4mm}
\includegraphics[height = 3.28cm, width = 3.6cm]{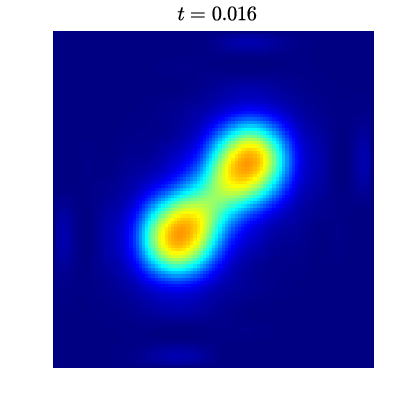}\hspace{-4mm}
\includegraphics[height = 3.28cm, width = 3.6cm]{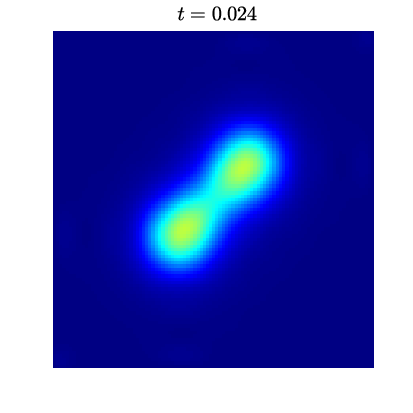}\hspace{-4mm}
\includegraphics[height = 3.28cm, width = 3.6cm]{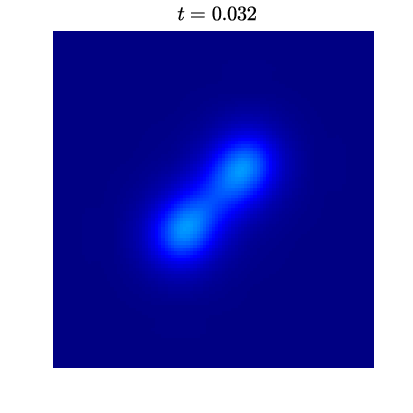}}\vspace{-2mm}
\centerline{
\hspace{3.3mm}
\includegraphics[height = 3.28cm, width = 3.6cm]{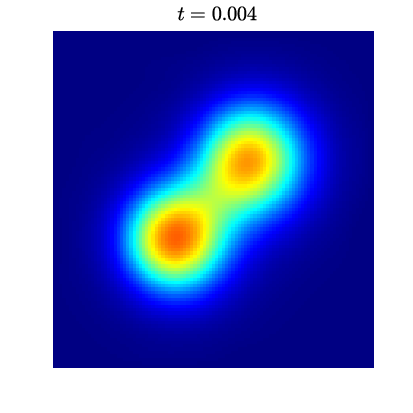}\hspace{-4mm}
\includegraphics[height = 3.28cm, width = 3.6cm]{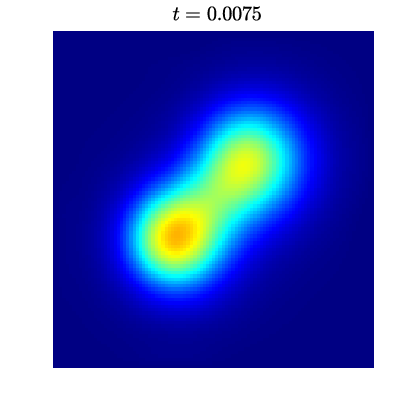}\hspace{-4mm}
\includegraphics[height = 3.28cm, width = 3.6cm]{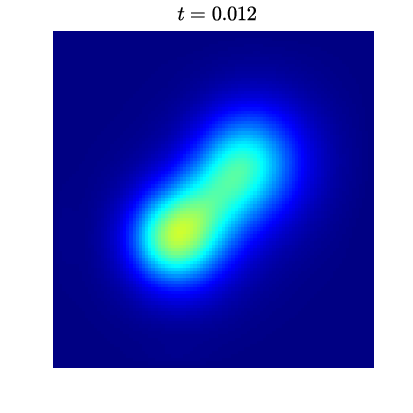}\hspace{-3mm}
\includegraphics[height = 3.38cm, width = 4.1cm]{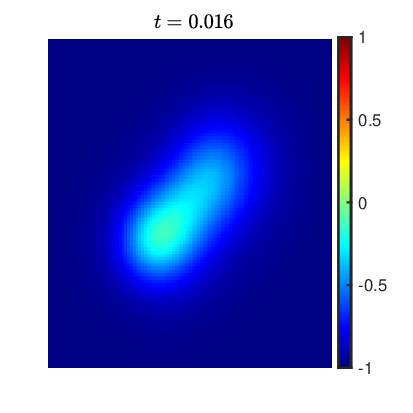}}
\caption{Time evolution of two bubbles in the Allen--Cahn equations (\ref{AC-1})--(\ref{u0-two_bubble}) with $\ap(\bx) \equiv 2$ (top row), $\ap(\bx) \equiv 1.5$ (middle row), $\ap(\bx) = 1.5+x/2$ (last row).}\label{Figure6-3-1}
\end{figure}

Figure \ref{Figure6-3-1} shows the time evolution of two bubbles in the Allen--Cahn equations with different $\ap(\bx)$. 
In the classical case,  two bubbles first coalesce into one (see Figure \ref{Figure6-3-1} top row for $t = 0.0075$), and then it is eventually absorbed by the fluid. 
While in the fractional case with $\ap(\bx) \equiv 1.5$,  the anomalous diffusion slows down the evolution of two bubbles \cite{Duo2019, Duo-FDM2019}.  
No complete coalescence is observed, and two bubbles diffuse separately and vanish  after a longer time. 
In contrast to symmetric evolution,  two bubbles behave differently in heterogeneous fluids. 
For  $\ap(\bx) = 1.5+x/2$,  the bubble at the right hand side is absorbed first. 
Generally, the larger the value of $\ap(\bx)$, the faster the bubble diminishes. 

\section{Conclusions}
\label{section6}
\setcounter{equation}{0}

We carried out an extensive study on the variable-order Laplacian $(-\Dt)^{\ap(\bx)/2}$ for  $0 < \ap(\bx) \le 2$.
The variable-order (fractional) Laplacian plays a significant role in modeling and studying heterogeneous systems. 
Surprisingly, it is still difficult to find adequate analytical and numerical studies on this operator in the current literature. 
We discussed two definitions of the variable-order fractional Laplacian  including the pseudo-differential definition in (\ref{pseudo}) and hypersingular integral definition in (\ref{integralFL}). 
It showed that in the special case of constant order, this operator reduces  to the well-known fractional Laplacian associated with the L\'evy process. 
We then  presented a class of hypergeometric functions whose variable-order Laplacian can be analytically expressed. 
These results further lead to the variable-order Laplacian of Gaussian functions, generalized inverse multiquadric functions,  and Bessel-type functions, all of which are well-known candidates for the positive definite radial basis functions. 
Therefore, these analytical results are important building blocks in developing our meshless methods, and they can also serve as benchmarks in designing and testing number methods for variable-order Laplacian.

We proposed a class of meshfree methods to solve problems with the variable-order Laplacian   $(-\Dt)^{\ap(\bx)/2}$. 
To the best of our knowledge, these are the  first numerical methods for the variable-order  (fractional) Laplacian. 
Our meshfree methods, integrating the advantages of both pseudo-differential representation and hypersingular integral form of the variable-order fractional Laplacian,  can solve heterogeneous problems  in a seamless manner. 
Moreover, utilizing the analytical results of RBFs, our methods bypass numerical approximation of the hypersingular integral of  fractional Laplacian and thus avoid large computational cost in evaluating fractional derivatives of RBFs. 
Our methods are simple and flexible of domain geometry, and their computer implementation remains the same for any dimension $d \ge 1$. 

Numerical studies in approximating Laplacian operators and solving Poisson problems showed that our method can achieve higher accuracy with fewer points in comparison to the finite difference method. 
This is very important in the study of fractional derivatives with variable order, as their heterogeneity greatly increases the storage and computational cost in numerical simulations. 
We then applied the proposed method to  study solution behaviors of variable-order fractional PDEs arising in different applications \cite{Zhu2014, Xue2018, Chen2014b, Lenzi2016, Delia2021}.
The transition of waves between classical and fractional media showed that the wave properties could be significantly changed in heterogeneous media. 
While the coexistence of anomalous and normal diffusion leads to an asymmetric solution propagation. 
These studies  provided insights for further understanding and applications of variable-order fractional derivatives. 
The study of variable-order fractional Laplacian still remains limited. 
In the future, we will carry out more mathematical and numerical research to further understand this heterogeneous operator.

\bibliographystyle{plain}

\begin{thebibliography}{99}

\bibitem{Acosta2017}
G.~Acosta and J.~P. Borthagaray.
\newblock A fractional {L}aplace equation: Regularity of solutions and finite
  element approximations.
\newblock {\em SIAM J. Numer. Anal.}, 55(2):472--495, 2017.

\bibitem{Ainsworth2018B}
M.~Ainsworth and C.~Glusa.
\newblock Towards an efficient finite element method for the integral
  fractional {L}aplacian on polygonal domains.
\newblock In {\em Contemporary computational mathematics---a celebration of the
 80th birthday of {I}an {S}loan. {V}ol. 1}:17--57. Springer, Cham,
  2018.

\bibitem{Antoine2021}
X.~Antoine, E.~Lorin, and Y.~Zhang.
\newblock Derivation and analysis of computational methods for fractional 
Laplacian equations with absorbing layers.
\newblock {\em Numer. Algorithms}, 87:409--444, 2021.

\bibitem{Baeumer2010}
B.~Baeumer and M.~M. Meerschaert.
\newblock Tempered stable {L}\'{e}vy motion and transient super-diffusion.
\newblock {\em J. Comput. Appl. Math.}, 233(10):2438--2448, 2010.

\bibitem{Bass1988}
R.~F. Bass.
\newblock Uniqueness in law for pure jump {M}arkov processes.
\newblock {\em Probab. Theory Rel.}, 79(2):271--287, 1988.


\bibitem{Bonito2019}
A.~Bonito, W.~Lei, and J.~E. Pasciak.
\newblock Numerical approximation of the integral fractional {L}aplacian.
\newblock {\em Numer. Math.}, 142(2):235--278, 2019.

\bibitem{Burkardt2021}
J.~Burkardt, Y.~Wu, and Y.~Zhang.
\newblock A unified meshfree pseudospectral method for solving both classical
  and fractional {PDE}s.
\newblock {\em SIAM J. Sci. Comput.}, 43(2):A1389--A1411, 2021.



\bibitem{Chen2014b}
H.~Chen, H.~Zhou, and S.~Qu.
\newblock Low rank approximation for time domain viscoacoustic wave equation
  with spatially varying order fractional laplacians.
\newblock {\em 84th Annual International Meeting}, SEG:3400--3405, 2014.

\bibitem{Chen2020}
X.~Chen, Chen Z.-Q., and J.~Wang.
\newblock Heat kernel for nonlocal operators with variable-order.
\newblock {\em Stoch. Proc. Appl.}, 130(6):3574--3647, 2020.

\bibitem{CRUZ-URIBE2013}
D.~V. Cruz-Uribe and A.~Fiorenza.
\newblock {\em Variable {L}ebesgue Spaces}.
\newblock Applied and Numerical Harmonic Analysis. Birkh\"{a}user/Springer,
  Heidelberg, 2013. 

\bibitem{Delia2021}
M.~D'Elia and C.~Glusa.
\newblock A fractional model for anomalous diffusion with increased variability: Analysis, algorithms and applications to interface problems. 
\newblock {\em Numer. Methods Partial Differ. Eq.}, 1--20, 2021.

\bibitem{DIENING2017}
L.~Diening, P.~Harjulehto, P.~H\"{a}st\"{o}, and M.~R\o{c}irc{u}\v{z}i\v{c}ka.
\newblock {\em Lebesgue and {S}obolev Spaces with Variable Exponents}, volume
  2017 of {\em Lecture Notes in Mathematics}.
\newblock Springer, Heidelberg, 2011.

\bibitem{Du2022}
R.~Du, Z.~Sun, and H.~Wang, 
\newblock Temporal second-order finite difference schemes for variable-order time-fractional wave equations. 
\newblock {\em SIAM J. Num. Anal.}, 60(1):104--132, 2022.

\bibitem{Dubrulle1998}
B.~Dubrulle and J.-P. Laval.
\newblock Truncated L\'evy laws and 2d turbulence.
\newblock {\em Phys. J. B}, (4):143--146, 1998.

\bibitem{Duo2018}
S.~Duo, H.~W. van Wyk, and Y.~Zhang.
\newblock A novel and accurate finite difference method for the fractional
  {L}aplacian and the fractional {P}oisson problem.
\newblock {\em J. Comput. Phys.}, 355:233--252, 2018.

\bibitem{Duo2015}
S.~Duo and Y.~Zhang.
\newblock Computing the ground and first excited states of the fractional {S}chr\"odinger equation in an infinite potential well.
\newblock {\em Commun. Comput. Phys.}, 18(2):321--350, 2015.

\bibitem{Duo2016}
S.~Duo and Y.~Zhang.
\newblock Mass-conservative {F}ourier spectral methods for solving the
  fractional nonlinear {S}chr\"{o}dinger equation.
\newblock {\em Comput. Math. with Appl.}, 71(11):2257--2271, 2016.

\bibitem{Duo-FDM2019}
S.~Duo and Y.~Zhang.
\newblock Accurate numerical methods for two and three dimensional integral
  fractional {L}aplacian with applications.
\newblock {\em Comput. Methods. Appl. Mech. Eng.}, 355:639--662, 2019.

\bibitem{Duo2019}
S.~Duo and H.~Wang.
\newblock A fractional phase-field model using an infinitesimal generator 
of $\alpha$ stable L\'evy process.
\newblock {\em J. Comput. Phys.}, 384:253--269, 2019.

\bibitem{Duo-TFL2019}
S.~Duo and Y.~Zhang.
\newblock Numerical approximations for the tempered fractional {L}aplacian:
  Error analysis and applications.
\newblock {\em J. Sci. Comput.}, 81(1):569--593, 2019.

\bibitem{Dwivedi2021}
K.~D. Dwivedi, Rajeev, S.~Das, and J.~F. Gomez-Aguilar.
\newblock Finite difference/collocation method to solve multi term
  variable-order fractional reaction–advection–diffusion equation in
  heterogeneous medium.
\newblock {\em Numer. Methods Partial Differ. Eq.}, 37(3):2031--2045, 2021.

\bibitem{Dyda2011}
B.~Dyda.
\newblock Fractional {H}ardy inequality with a remainder term.
\newblock {\em Colloq. Math.}, 122(1):59--67, 2011.

\bibitem{Dyda2012}
B.~Dyda.
\newblock Fractional calculus for power functions and eigenvalues of the
  fractional {L}aplacian.
\newblock {\em Fract. Calc. Appl. Anal.}, 15(4):536--555, 2012.

\bibitem{Dyda2017}
B.~Dyda, A.~Kuznetsov, and M.~Kwa\'{s}nicki.
\newblock Fractional {L}aplace operator and {M}eijer {G}-function.
\newblock {\em Constr. Approx.}, 45(3):427--448, 2017.

\bibitem{Flyer2006}
N.~Flyer.
\newblock Exact polynomial reproduction for oscillatory radial basis functions
  on infinite lattices.
\newblock {\em Comput. Math. Appl.}, 51(8):1199--1208, 2006.

\bibitem{Fornberg2002}
B.~Fornberg, T.~A. Driscoll, G.~Wright, and R.~Charles.
\newblock Observations on the behavior of radial basis function approximations
  near boundaries.
\newblock {\em Comput. Math. Appl.}, 43(3-5):473--490, 2002.

\bibitem{Fornberg2006}
B.~Fornberg, E.~Larsson, and G.~Wright.
\newblock A new class of oscillatory radial basis functions.
\newblock {\em Comput. Math. Appl.}, 51(8):1209--1222, 2006.

\bibitem{Fornberg2017}
B.~Fornberg and J.~Zuev.
\newblock The {R}unge phenomenon and spatially variable shape parameters in
  {RBF} interpolation.
\newblock {\em Comput. Math. Appl.}, 54(3):379--398, 2007.

\bibitem{GRADSHTEYN2007}
I.~S. Gradshteyn and I.~M. Ryzhik.
\newblock {\em Table of Integrals, Series, and Products}.
\newblock Academic Press, Amsterdam, seventh edition, 2007.

\bibitem{Hao2021}
Z.~Hao, Z.~Zhang, and R.~Du, 
\newblock Finite centered difference scheme for high-dimensional integral fractional Laplacian, 
\newblock{\em J. Comput. Phys.}, 424:109851, 2021.

\bibitem{Hormander1965}
L.~H{\"o}rmander.
\newblock Pseudo-differential operators.
\newblock {\em Comm. Pure Appl. Math.}, 18:501--517, 1965.

\bibitem{Javanainen2013}
M.~Javanainen, H.~Hammar{\'e}n, L.~Monticelli, J.~Jeon, M.~S. Miettinen,
  H.~Martinez-Seara, R.~Metzler, and I.~Vattulainen.
\newblock Anomalous and normal diffusion of proteins and lipids in crowded
  lipid membranes.
\newblock {\em Faraday Discuss.}, 161:397--417, 2013.

\bibitem{Kansa90I}
E.~J. Kansa.
\newblock Multiquadrics -- A scattered data approximation scheme with
  applications to computational fluid-dynamics. {I}. {S}urface approximations
  and partial derivative estimates.
\newblock {\em Comput. Math. Appl.}, 19(8--9):127--145, 1990.

\bibitem{Kansa90II}
E.~J. Kansa.
\newblock Multiquadrics -- A scattered data approximation scheme with
  applications to computational fluid-dynamics. {II}. {S}olutions to parabolic,
  hyperbolic and elliptic partial differential equations.
\newblock {\em Comput. Math. Appl.}, 19(8--9):147--161, 1990.

\bibitem{Kikuchi1997}
K.~Kikuchi and A.~Negoro.
\newblock On {M}arkov process generated by pseudodifferential operator of
  variable order.
\newblock {\em Osaka J. Math.}, 34(2):319--335, 1997.


\bibitem{Kirkpatrick2016}
K.~Kirkpatrick and Y.~Zhang.
\newblock Fractional {S}chr\"{o}dinger dynamics and decoherence.
\newblock {\em Physica D}, 332:41--54, 2016.

\bibitem{Kohn1965}
J.~J. Kohn and L.~Nirenberg.
\newblock An algebra of pseudo-differential operators.
\newblock {\em Comm. Pure Appl. Math.}, 18(1-2):269--305, 1965.

\bibitem{Kuhn2017}
F.~K\"{u}hn.
\newblock {\em L\'evy Type Processes: Moments, Construction and Heat Kernel Estimates}.
\newblock Springer Lecture Notes in Mathematics, vol. 2187. Springer, Berlin, 2017.

\bibitem{Kuhn2020}
F.~K\"{u}hn.
\newblock Schauder estimates for 
Poisson equations associated with non-local
  feller generators.
\newblock {\em J. Theor. Probab.}, 17, 2020.



\bibitem{Larsson2005}
E.~Larsson and B.~Fornberg.
\newblock Theoretical and computational aspects of multivariate interpolation
  with increasingly flat radial basis functions.
\newblock {\em Comput. Math. Appl.}, 49(1):103--130, 2005.

\bibitem{Lenzi2016}
E.~K. Lenzi, H.~V. Ribeiro, A.~A. Tateishi, R.~S. Zola, and L.~R. Evangelista.
\newblock Anomalous diffusion and transport in heterogeneous systems separated
  by a membrane.
\newblock {\em Proc. A.}, 472(2195):20160502, 2016.

\bibitem{Leopold1999}
H.~Leopold.
\newblock Embedding of function spaces of variable order of differentiation in
  function spaces of variable order of integration.
\newblock {\em Czechoslovak Math. J.}, 49(124)(3):633--644, 1999.

\bibitem{Lin2009}
R.~Lin, F.~Liu, V.~Anh, and I.~Turner.
\newblock Stability and convergence of a new explicit finite-difference
  approximation for the variable-order nonlinear fractional diffusion equation.
\newblock {\em Appl. Math. Comput.}, 212(2):435--445, 2009.

\bibitem{Lorenzo2002}
C.~F. Lorenzo and T.~T. Hartley.
\newblock Variable order and distributed order fractional operators.
\newblock {\em Nonlinear Dynam.}, 29(1--4):57--98, 2002.

\bibitem{Luo2019}
D.~Luo and J.~Wang.
\newblock Coupling by reflection and {H}\"{o}lder regularity for non-local
  operators of variable order.
\newblock {\em Trans. Amer. Math. Soc.}, 371(1):431--459, 2019.

\bibitem{Meerschaert2006}
M.~M. Meerschaert and C.~Tadjeran.
\newblock Finite difference approximations for two-sided space-fractional
  partial differential equations.
\newblock {\em Appl. Numer. Math.}, 56(1):80--90, 2006.

\bibitem{Meerschaert2008}
M.~M. Meerschaert, Y.~Zhang, and B.~Baeumer.
\newblock Tempered anomalous diffusion in heterogeneous systems.
\newblock {\em Geophys. Res. Lett.}, 35(17):L17403, 2008.


\bibitem{Pang2015}
G.~Pang, W.~Chen, and Z.~Fu.
\newblock Space-fractional advection-dispersion equations by the {K}ansa
  method.
\newblock {\em J. Comput. Phys.}, 293:280--296, 2015.

\bibitem{Piret2013}
C.~Piret and E.~Hanert.
\newblock A radial basis functions method for fractional diffusion equations.
\newblock {\em J. Comput. Phys.}, 238:71--81, 2013.

\bibitem{PRUDNIKOV1990}
A.~P. Prudnikov, Yu.~A. Brychkov, and O.~I. Marichev.
\newblock {\em Integrals and series. {V}ol. 3}.
\newblock Gordon and Breach Science Publishers, New York, 1990.

\bibitem{Rafeiro2016}
H.~Rafeiro and S.~G. Samko.
\newblock Fractional integrals and derivatives: Mapping properties.
\newblock {\em Fract. Calc. Appl. Anal.}, 19(3):580--607, 2016.

\bibitem{Rosenfeld2019}
J.~A. Rosenfeld, S.~A. Rosenfeld, and W.~E. Dixon.
\newblock A mesh-free pseudospectral approach to estimating the fractional
  {L}aplacian via radial basis functions.
\newblock {\em J. Comput. Phys.}, 390:306--322, 2019.

\bibitem{Samko2013}
S.~G. Samko.
\newblock Fractional integration and differentiation of variable order: An
  overview.
\newblock {\em Nonlinear Dynam.}, 71(4):653--662, 2013.

\bibitem{SAMKO1993b}
S.~G. Samko, A.~A. Kilbas, and O.~I. Marichev.
\newblock {\em Fractional Integrals and Derivatives}.
\newblock Gordon and Breach Science Publishers, Yverdon, 1993.

\bibitem{Samko1993}
S.~G. Samko and B.~Ross.
\newblock Integration and differentiation to a variable fractional order.
\newblock {\em Integral Transform. Spec. Funct.}, 1(4):277--300, 1993.

\bibitem{Sarra2009}
S.~A. Sarra and E.~J. Kansa.
\newblock Multiquadric radial basis function approximation methods for the
  numerical solution of partial differential equations.
\newblock {\em Adv. in Comput. Mech.}, 2, 2009.


\bibitem{Song2016}
F.~Song, C.~Xu, and G.~E. Karniadakis.
\newblock A fractional phase-field model for two-phase flows with tunable
  sharpness:  Algorithms and simulations.
\newblock {\em Comput. Methods Appl. Mech. Engrg.}, 305:376--404, 2016.

\bibitem{Sun2019}
H.~Sun, A.~Chang, Y.~Zhang, and W.~Chen.
\newblock A review on variable-order fractional differential equations:
  Mathematical foundations, physical models, numerical methods and
  applications.
\newblock {\em Fract. Calc. Appl. Anal.}, 22(1):27--59, 2019.

\bibitem{Tang2020}
T.~Tang, L.-L.~Wang, H.~Yuan, and T.~Zhou.
\newblock Rational spectral methods for PDEs involving fractional Laplacian 
in unbounded domains.
\newblock {\em SIAM J. Sci. Comput.}, 42(2):A585--A611, 2020.

\bibitem{Tsuchiya1992}
M.~Tsuchiya.
\newblock L\'{e}vy measure with generalized polar decomposition and the
  associated {SDE} with jumps.
\newblock {\em Stochastics Stochastics Rep.}, 38(2):95--117, 1992.

\bibitem{Wu0020}
Y.~Wu and Y.~Zhang.
\newblock A universal solution scheme for fractional and classical {PDE}s.
\newblock {\em arXiv:2102.00113}, 2020.

\bibitem{Xiang2019}
M.~Xiang, B.~Zhang, and D.~Yang.
\newblock Multiplicity results for variable-order fractional {L}aplacian
  equations with variable growth.
\newblock {\em Nonlinear Anal.}, 178:190--204, 2019.

\bibitem{Xue2018}
Z.~Xue, H.~Baek, H.~Zhang, Y.~Zhao, T.~Zhu, and S.~Fomel.
\newblock Solving fractional {L}aplacian viscoelastic wave equations using domain
  decomposition.
\newblock {\em 88th Annual International Meeting}, SEG:3943--3947, 2018.

\bibitem{Zhang2012}
Y.~Zhang, M.~Meerschaert, and A.~Packman.
\newblock Linking fluvial bed sediment transport across scales.
\newblock {\em Geophys. Res. Lett.}, 39:L20404, 2012.

\bibitem{Zhao2015}
X.~Zhao, Z.~Sun, and G.~E. Karniadakis.
\newblock Second-order approximations for variable order fractional
  derivatives: Algorithms and applications.
\newblock {\em J. Comput. Phys.}, 293:184--200, 2015.

\bibitem{Zhu2014}
T.~Zhu and J.~M. Harris.
\newblock Modeling acoustic wave propagation in heterogeneous attenuating media
  using decoupled fractional {L}aplacians.
\newblock {\em Geophysics}, 79(3):T105--T116, 2014.

\end{thebibliography}

%
%
\end{document}